\documentclass[leqno,a4paper,12pt]{article}
\usepackage{amsmath,amsthm}

\newtheorem{Thm}{\indent Theorem}[section]
\newtheorem{Prop}[Thm]{\indent Proposition}
\newtheorem{Lem}[Thm]{\indent Lemma}
\newtheorem{Cor}[Thm]{\indent Corollary}

\theoremstyle{definition}
\newtheorem{Def}[Thm]{\indent Definition}
\newtheorem{Rem}[Thm]{\indent Remark}

\newtheorem{Ass}[Thm]{\indent Assumption}
\def\qed{{\hskip0pt\unskip\unskip\nobreak\hfil\penalty50
          \hskip1em\hbox{}\nobreak\hfil
          {\bf q.e.d.}%
          \parfillskip=0pt\finalhyphendemerits=0
          \par}\medskip}

\newenvironment{Proof}
               {{\it Proof.}\quad}
               {\qed}

\newenvironment{Proofof}[1]
               {{\it Proof of #1.}\quad}
               {\qed}

\newcommand{\Prime}{\kern3\fontdimen1\font$'$\kern-7\fontdimen1\font}

\long\def\forget#1{}

\long\def\beginSIDEREMARK#1\endSIDEREMARK
    {{\par\bigskip\advance\leftskip by 2cm
                  \advance\rightskip by -2cm\noindent
      {\bf Our own side remark:} #1
      \par\bigskip\noindent}}
\long\def\beginFORGET#1\endFORGET{#1}
\long\def\beginFORGET#1\endFORGET{}

\def\?{\ ???\ \immediate\write16{}%
\immediate\write16{Warning: There was still a question mark . . . }%
\immediate\write16{}}

\usepackage{amsmath}
\usepackage{exscale}
\usepackage{amssymb}

\newcommand{\BA}{{\mathbb{A}}}

\newcommand{\BC}{{\mathbb{C}}}
\newcommand{\BD}{{\mathbb{D}}}

\newcommand{\BG}{{\mathbb{G}}}

\newcommand{\BQ}{{\mathbb{Q}}}
\newcommand{\BR}{{\mathbb{R}}}

\newcommand{\BZ}{{\mathbb{Z}}}

\newcommand{\Fc}{{\mathfrak{c}}}

\newcommand{\Fn}{{\mathfrak{n}}}

\newcommand{\FH}{{\mathfrak{H}}}

\newcommand{\FS}{{\mathfrak{S}}}

\newcommand{\FX}{{\mathfrak{X}}}
\newcommand{\FY}{{\mathfrak{Y}}}
\newcommand{\FZ}{{\mathfrak{Z}}}

\newcommand{\CC}{{\cal C}}
\newcommand{\CD}{{\cal D}}

\newcommand{\CH}{{\cal H}}

\newcommand{\CK}{{\cal K}}
\newcommand{\CL}{{\cal L}}
\newcommand{\CM}{{\cal M}}

\newcommand{\CO}{{\cal O}}

\newcommand{\CV}{{\cal V}}
\newcommand{\CW}{{\cal W}}

\newcommand{\ch}{\mathop{\rm CH}\nolimits}
\newcommand{\Spec}{\mathop{{\bf Spec}}\nolimits}

\newcommand{\imm}{\mathop{{\rm im}}\nolimits}

\newcommand{\End}{\mathop{\rm End}\nolimits}

\newcommand{\GL}{\mathop{\rm GL}\nolimits}
\newcommand{\Gm}{\mathop{\BG_m}\nolimits}
\newcommand{\GKm}{\mathop{\BG_{m,L}}\nolimits}
\newcommand{\Gr}{\mathop{\rm Gr}\nolimits}

\newcommand{\Hom}{\mathop{\rm Hom}\nolimits}

\newcommand{\Res}{\mathop{\rm Res}\nolimits}
\newcommand{\Sym}{\mathop{\rm Sym}\nolimits}

\newcommand{\sign}{\mathop{{\rm sign}}\nolimits}
\newcommand{\loccit}{[loc.$\;$cit.]}

\def\halb{\frac{1}{2}}

\newbox\mybox
\def\arrover#1{\mathrel{
       \setbox\mybox=\hbox spread 1.4em{\hfil$\scriptstyle#1$\hfil}
       \vbox{\offinterlineskip\copy\mybox
             \hbox to\wd\mybox{\rightarrowfill}}}}
\def\larrover#1{\mathrel{
       \setbox\mybox=\hbox spread 1.4em{\hfil$\scriptstyle#1$\hfil}
       \vbox{\offinterlineskip\copy\mybox
             \hbox to\wd\mybox{\leftarrowfill}}}}

\def\ontoover#1{\mathrel{
       \setbox\mybox=\hbox spread 1.4em{\hfil$\scriptstyle#1$\hfil}
       \vbox{\offinterlineskip\copy\mybox
             \hbox to\wd\mybox{\rightarrowfill\hskip-2.8mm
                               $\rightarrow$}}}}
\def\leftontoover#1{\mathrel{
       \setbox\mybox=\hbox spread 1.4em{\hfil$\scriptstyle#1$\hfil}
       \vbox{\offinterlineskip\copy\mybox
             \hbox to\wd\mybox{$\leftarrow$\hskip-2.8mm
                               \leftarrowfill}}}}
\def\longto{\longrightarrow}
\def\into{\hookrightarrow}

\def\isoto{\arrover{\sim}}

\def\longinto{\lhook\joinrel\longrightarrow}

\usepackage[curve,matrix,arrow,cmtip]{xy}
\NoComputerModernTips

\def\myxymessage{\def\messagetext
   {Here an xy-pic diagram was omitted to speed up compilation . . . }
   \immediate\write16{\messagetext}
   \hbox{\bf \messagetext}}
\def\filxymatrix#1{\myxymessage}
\def\filxyarray#1{\myxymessage}
\newdir^{ (}{{}*!/-3pt/\dir^{(}}
\newdir_{ (}{{}*!/-3pt/\dir_{(}}
\newdir^{ )}{{}*!/+3pt/\dir^{)}}
\newdir_{ )}{{}*!/+3pt/\dir_{)}}

\def\rscript#1{\hbox to 0pt{$\scriptstyle#1$\hss}}

\let\oldbullet\bullet
\def\bullet{{\mathchoice{\oldbullet}%
                        {\oldbullet}%
                        {\scriptscriptstyle\oldbullet}%
                        {\oldbullet}}}

\newcommand{\bA}{\mathop{\overline{A}}\nolimits}
\newcommand{\bbA}{\mathop{\overline{\overline{A}}}\nolimits}
\newcommand{\bS}{\mathop{\overline{S}}\nolimits}

\newcommand{\bX}{\mathop{\overline{X}}\nolimits}

\newcommand{\ur}{\mathop{\underline{r}}\nolimits}

\newcommand{\CHeffM}{\mathop{CHM^{eff}(k)}\nolimits}

\newcommand{\CHM}{\mathop{CHM(k)}\nolimits}
\newcommand{\CHSM}{\mathop{CHM(S)}\nolimits}

\newcommand{\CHeffQM}{\mathop{\CHeffM_F}\nolimits}
\newcommand{\CHQM}{\mathop{\CHM_F}\nolimits}
\newcommand{\CHSQM}{\mathop{CHM(S)_F}\nolimits}
\newcommand{\DeffgM}{\mathop{DM^{eff}_{gm}(k)}\nolimits}
\newcommand{\DeffQgM}{\mathop{\DeffgM_F}\nolimits}

\newcommand{\DgM}{\mathop{DM_{gm}(k)}\nolimits}
\newcommand{\DQgM}{\mathop{\DgM_F}\nolimits}

\newcommand{\QDM}{\mathop{M \! D(k)_F}\nolimits}
\newcommand{\QATM}{\mathop{M \! AT(k)_F}\nolimits}
\newcommand{\DQATM}{\mathop{D \! \QATM}\nolimits}

\newcommand{\ReL}{\mathop{{\rm Res}_{L/\BQ}}\nolimits}
\newcommand{\Mgm}{\mathop{M_{gm}}\nolimits}

\newcommand{\Mcgm}{\mathop{M_{gm}^c}\nolimits}

\newcommand{\dMgm}{\mathop{\partial M_{gm}}\nolimits}
\begin{document}

%%%%%%%%%%%%%%%%%%%%%%%%%%%%%%%%%%%%%%%%%%%%%%%%%%%%%%%%%%%%%%%%%%%%%%%
%
%  formatting

\hfuzz=3pt
\overfullrule=10pt                   % erzeugt schwarze Fehlerbalken

% The displayskip values were changed because \LaTeX does not react
% correctly to a \leqno: it should then use big skips, but doesn't.

\setlength{\abovedisplayskip}{6.0pt plus 3.0pt}
                               % preset 10.0pt plus 2.0pt minus 5.0pt
\setlength{\belowdisplayskip}{6.0pt plus 3.0pt}
                               % preset 10.0pt plus 2.0pt minus 5.0pt
\setlength{\abovedisplayshortskip}{6.0pt plus 3.0pt}
                               % preset 0.0pt plus 3.0pt
\setlength{\belowdisplayshortskip}{6.0pt plus 3.0pt}
                               % preset 6.0pt plus 3.0pt minus 3.0pt

\setlength{\baselineskip}{13.0pt}
                               % preset 12.0pt
\setlength{\lineskip}{0.0pt}
                               % preset 1.0pt
\setlength{\lineskiplimit}{0.0pt}
                               % preset 0.0pt

%%%%%%%%%%%%%%%%%%%%%%%%%%%%%%%%%%%%%%%%%%%%%%%%%%%%%%%%%%%%%%%%%%%%%%%
%
%  Title Page
%
%%%%%%%%%%%%%%%%%%%%%%%%%%%%%%%%%%%%%%%%%%%%%%%%%%%%%%%%%%%%%%%%%%%%%%%

\title{On the interior motive of certain Shimura varieties:
the case of Hilbert--Blumenthal varieties
\forget{
\footnotemark
\footnotetext{To appear in ....}
}
}
\author{\footnotesize by\\ \\
\mbox{\hskip-2cm
\begin{minipage}{6cm} \begin{center} \begin{tabular}{c}
J\"org Wildeshaus \footnote{
Partially supported by the \emph{Agence Nationale de la
Recherche}, project no.\ ANR-07-BLAN-0142 ``M\'ethodes \`a la
Voevodsky, motifs mixtes et G\'eom\'etrie d'Arakelov''. }\\[0.2cm]
\footnotesize LAGA\\[-3pt]
\footnotesize UMR~7539\\[-3pt]
\footnotesize Institut Galil\'ee\\[-3pt]
\footnotesize Universit\'e Paris 13\\[-3pt]
\footnotesize Avenue Jean-Baptiste Cl\'ement\\[-3pt]
\footnotesize F-93430 Villetaneuse\\[-3pt]
\footnotesize France\\
{\footnotesize \tt wildesh@math.univ-paris13.fr}
\end{tabular} \end{center} \end{minipage}
\hskip-2cm}
\\[2.5cm]
\forget{
{\bf Preliminary version --- not for distribution}\\[1cm]
}
}
% In the final version we might want to fix the date:
\date{March 18, 2011}
\maketitle
%\quad \\[-1.7cm]
\begin{abstract}
\noindent
The purpose of this article is
to construct a Hecke-equivariant Chow motive whose realizations equal
interior (or intersection) cohomology of Hilbert--Blumenthal 
varieties with non-constant algebraic coefficients.  \\


%\noindent
%{\bf R\'esum\'e~:} RESUME.\\
\end{abstract}

%\vfill

\bigskip
\bigskip
\bigskip

\noindent {\footnotesize Math.\ Subj.\ Class.\ (2000) numbers: 14G35
(11F32, 11F41, 14C25, 14F25, 19E15, 19F27)
}

\eject
\tableofcontents

\bigskip
%\vspace*{0.5cm}

%\newpage
%\include{Intro}

%%%%%%%%%%%%%%%%%%%%%%%%%%%%%%%%%%%%%%%%%%%%%%%%%%%%%%%%%%%%%%%%%%%%%%%
%
%  Introduction
%
%%%%%%%%%%%%%%%%%%%%%%%%%%%%%%%%%%%%%%%%%%%%%%%%%%%%%%%%%%%%%%%%%%%%%%%

\setcounter{section}{0}
\section{Introduction}
\label{Intro}

%%%%%%%%%%%%%%%%%%%%%%%%%%%%%%%%

%%%%%%%%%%%%%%%%%%%%%%%%%%%%%%%%

The purpose of this paper is the construction and the analysis 
of the \emph{interior motive} of Kuga--Sato families over
Hilbert--Blumenthal varieties associated to a 
fixed totally real number field $L$. \\

In order to motivate the problem, let us start by discussing 
the case of classical modular curves, i.e., of 
Hilbert--Blumenthal varieties associated to $L = \BQ \, $.
Fix two integers $n \ge 3$ and $r \ge 1$, let $S$
denote the modular curve parametrizing elliptic curves
with level $n$ structure, and $S \into \bS$ its smooth
compactification.
Write $A \to S$ for the universal elliptic curve, and 
$\bA \to \bS$ for the universal generalized elliptic
curve. Thus, $\bA$ is smooth and proper over $\BQ \, $. The $r$-fold
fibre product 
$\bA^r := \bA \times_{\bS} \times \ldots \times_{\bS} \bA$ 
of $\bA$ over $\bS$ is singular
for $r \ge 2$,
and can be desingularized canonically \cite{D}. 
Denote by $\bbA^r$ this desingularization.
The symmetric group $\FS_r$ acts on $\bA^r$ by permutations,
the $r$-th power of the group $\BZ / n \BZ$ by translations, and the
$r$-th power of the
group $\mu_2$ by inversion in the fibres. Altogether, 
this gives a canonical action of the semi-direct product
\[
\Gamma_r := \bigl( (\BZ / n \BZ)^2 \rtimes \mu_2 \bigr)^r \rtimes \FS_r
\]
by automorphisms on $\bA^r$. By the canonical nature of the desingularisation,
this extends to an action of $\Gamma_r$ by automorphisms on $\bbA^r$. \\

Let $\varepsilon: \Gamma_r \to \{ \pm 1 \}$
be the morphism which is trivial on $(\BZ / n \BZ)^{2r}$, is the
product map on $\mu_2^r$, and is the sign character on $\FS_r$.
Let $e$ denote the idempotent in the group ring
$\BQ[\Gamma_r]$ associated to $\varepsilon$. 
Following \cite{Sch}, one defines the Chow motive ${}^r_n \CW$ as the image 
\[
{}^r_n \CW := \Mgm \bigl( \bbA^r \bigr)^e 
\]
of the idempotent $e$ on the Chow motive $\Mgm \bigl( \bbA^r \bigr)$
of the smooth and proper $\BQ$-scheme $\bbA^r$. \\

By \cite{Sch}, 
the (Betti or $\ell$-adic) realizations of ${}^r_n \CW$
equal \emph{interior cohomo\-logy}, i.e., the image of the morphism
\[
H^n_c \bigl( A^r \bigr)^e \longto 
H^n \bigl( A^r \bigr)^e 
\]
linking (Betti or $\ell$-adic) cohomology with compact support to cohomology
without support of the $r$-fold
fibre product $A^r$ of $A$ over $S$. 
Note that at places $p$ of good reduction of $\bbA^r$
(i.e., $p$ does not divide $n$),
this result implies reduction properties for the Galois 
action on $\ell$-adic interior cohomology \cite[Thm.~1.2.4]{Sch}. \\

Using the Eichler--Shimura isomorphism, the Hodge structure or Galois module
associated to an elliptic normalized newform $f$ of 
level $n$ and weight $r + 2$ can be realized
in (the image of $e$ on) interior cohomology of $A^r$. This realization
coincides with the eigenspace of the action of the Hecke algebra
associated to the eigenvalues $a_p$ occurring in the Fourier expansion of $f$
\cite{D}. In order to lift this construction to motives,
one first extends the Hecke corres\-pondences to $\bbA^r$,
and checks that they commute with the action of $\Gamma_r$
\cite[Sect.~4.1]{Sch}.
Scholl then constructs the Grothendieck motive $M(f)$ associated to $f$
as the factor of the Grothendieck motive 
underlying the Chow motive ${}^r_n \CW$ 
corres\-ponding to the same eigenvalues $a_p$ \cite[Sect.~4.2]{Sch}. \\

Let us
agree on the principle that in order to geometrically explain
certain \emph{purity phenomena} in Hodge or Galois theory,
it is desirable to construct (Chow or Grothendieck) motives 
whose realizations equal the Hodge structure or Galois module
in question. \\

In the setting discussed above, purity concerns
interior cohomology of the smooth non-proper scheme $A^r$,
and more specifically, the weight occurring in the Hodge structure
or Galois module associated to the modular form $f$. 
The motives ${}^r_n \CW$
and $M(f)$ are constructed out of a smooth and proper scheme.
Given the definition of interior cohomology, it is natural
to look for this scheme among the smooth compactifications
of $A^r$. \\

Let us insist on one decisive additional point: 
among the smooth compactifications of $A^r$, there is one choice
(namely $\bbA^r$) allowing for extensions of the actions on $A^r$
of both the finite group $\Gamma_r$ and the Hecke correspondences. 
Unfortunately, this additional point is rather atypical:
smooth equivariant compactifications are not known 
(and maybe not reasonable to expect) to exist for Shimura varieties
of higher dimensions. Thus, it appears unrealistic to hope
for a na{\"{\i}}ve generalization of Scholl's approach
to such varieties.  \\

In the present paper, we shall develop an alternative strategy,
based on the formalism of \emph{weight structures} \cite{Bo}.
More precisely, we shall use the \emph{motivic weight structure}
of \loccit , with which
the category $\DeffgM$ of effective
\emph{geometrical motives} over a perfect field $k$ \cite{V} is equipped. \\

Since our strategy gives the same result as \cite{Sch}
for modular curves, let us continue to discuss that case,
and make the necessary translations. For a smooth scheme over $k$,
denote by $\Mgm (X)$ and $\Mcgm (X)$
the \emph{motive} of $X$ and its \emph{motive with compact support}, 
respectively \cite{V}. Effective Chow motives form a full sub-category
of $\DeffgM$; indeed they are identified with the category
of objects which are \emph{pure of weight zero} with respect to
the motivic weight structure \cite{Bo}.
Both functors $\Mgm$ and $\Mcgm$ agree on smooth and proper 
$k$-schemes $X$, and yield the Chow motive of $X$. 
In particular, we may consider ${}^r_n \CW$ as a geometrical
motive over $\BQ \, $. Thanks to the functoriality properties of
$\Mgm$ and $\Mcgm$, we get natural morphisms
\[
j^{r,*} : 
{}^r_n \CW = \Mgm \bigl( \bbA^r \bigr)^e \longto
\Mcgm \bigl( A^r \bigr)^e  
\]
and 
\[
j^r : 
\Mgm \bigl( A^r \bigr)^e \longto \Mgm \bigl( \bbA^r \bigr)^e = {}^r_n \CW
\]
induced by the open immersion $j^r$ of $A^r$ into $\bbA^r$. 
The analysis from \cite[Sect.~2 and 3]{Sch} of the geometric properties 
of the boundary of $\bbA^r$ can then be employed
to prove \cite[Sect.~3]{W3}
that $j^{r,*}$ and $j^r$ factor canonically through isomorphisms
\[
{}^r_n \CW \isoto \Gr_0 \Mcgm \bigl( A^r \bigr)^e 
\]
and 
\[
\Gr_0 \Mgm \bigl( A^r \bigr)^e \isoto {}^r_n \CW \; .
\] 
Here, the symbol $\Gr_0$ denotes the ``graded part of weight zero''
with respect to the motivic weight structure. 
It is defined on motives enjoying an additional property,
which therefore turns out to be essential for the present paper: the property of
\emph{avoiding weights} $-1$ and
$1$. Both $\Mcgm \bigl( A^r \bigr)^e$ and 
$\Mgm \bigl( A^r \bigr)^e$ satisfy this hypothesis, thanks
again to the analysis of the boundary of $\bbA^r$ \cite{Sch}. \\

To resume, we see that the canonical morphism
$\Mgm \bigl( A^r \bigr)^e \to 
\Mcgm \bigl( A^r \bigr)^e$ factors canonically through an isomorphism
\[ 
\Gr_0 \Mgm \bigl( A^r \bigr)^e \isoto
\Gr_0 \Mcgm \bigl( A^r \bigr)^e \; .
\]  
Our key observation is
that this latter statement does no longer necessitate the reference 
to a compactification of $A^r \, $! 
In fact, it is this statement that turns out to
generalize, with appropriate choices of $e$ (see below) to  
Hilbert--Blumenthal varieties of arbitrary dimension.
Since it concerns the ``graded part of weight zero'' of the
canonical morphism relating $\Mgm \bigl( A^r \bigr)^e$ to 
$\Mcgm \bigl( A^r \bigr)^e$, it necessitates control
of the weights present in a cone of that morphism. \\

Let us make this a little more precise.
Recall that for a smooth scheme $X$ over $k$,
the \emph{boundary motive} $\dMgm(X)$ of $X$ \cite{W1}
is a canonical choice of (the shift by [-1] of) such a cone; indeed, it
fits into a canonical exact triangle
\[
(\ast) \quad\quad
\dMgm(X) \longto \Mgm(X) \longto \Mcgm(X) \longto \dMgm(X)[1] \; .
\]
Given the specific choice of cone, the boundary motive satisfies
good functorial properties.
Assume that an idempotent endomorphism
$e$ of the exact triangle $(\ast)$ is given. We thus get a direct factor
\[
\dMgm(X)^e \longto \Mgm(X)^e \longto \Mcgm(X)^e \longto \dMgm(X)^e[1] \; .
\]  
In view of the above, 
let us assume the object $\dMgm(X)^e$ to avoid
weights $-1$ and $0$ with respect fo the motivic weight structure.
As verified in \cite{W3}, this assumption is not only
necessary, but also sufficient to construct the
Chow motives $\Gr_0 \Mgm(X)^e$ and $\Gr_0 \Mcgm(X)^e$,
and to show that they are canonically isomorphic. 
Given the nature of the realizations of $\Gr_0 \Mgm(X)^e$, it is natural to
call this object the $e$-part of the interior motive of $X$.  \\

In the context of Hilbert--Blumenthal varieties, it thus remains 
to verify the hypothesis on the absence of weights $-1$ and $0$ in 
$\dMgm(X)^e$, and technically speaking, this is what the present paper 
is about.
In Section~\ref{3}, we establish a criterion 
(Theorem~\ref{3Main}) allowing for that verification,
when the boundary motive is 
\emph{Artin--Tate}. Its simplified form (Corollary~\ref{3C})
states that the
absence of weights can be read off from the Hodge structure
or the Galois action on \emph{boundary cohomology}. \\

Section~\ref{4} contains the statements of
our main results, Theorems~\ref{4Main1} and \ref{4Main2}. Here,
$X$ is the $r$-th power of the universal Abelian scheme over a smooth
Hilbert--Blumenthal variety $S$ associated to a totally real number field
$L$ of degree $g$. The
idempotent $e$ cuts out the direct factor of the \emph{relative Chow motive}
of $X$, on which the action of $L$ is of type 
$(r_1,\ldots,r_g)$, for $r_1 + \ldots + r_g = r$. 
Theorem~\ref{4Main1} implies in particular
that in this context, the criterion from
Theorem~\ref{3Main} is satisfied
as soon as $r \ge 1$. Therefore, the interior motive $\Gr_0 \Mgm(X)^e$ exists. 
We list its principal pro\-perties, using the main results
from \cite[Sect.~4]{W3}.
First (Corollary~\ref{4E}), we get precise statements on the
weights occurring in the motive $\Mgm(X)^e$ and the motive with compact
support $\Mcgm(X)^e$. Second (Corollary~\ref{4Ea}),
the interior motive is Hecke-equivariant; it is here that the functorial
properties of the boundary motive turn out to be essential. Corollary~\ref{4Ea} appears
particularly interesting, given the problem of non-existence of
equivariant smooth compactifications of $X$ (for $g \ge 2$)
raised above.
\forget{
In this context, let us mention the main result of \cite{GHM1},
which implies the existence of relative Chow--K\"unneth projectors
$\Pi_S^i$, whose (relative) Betti realizations are isomorphic to the intersection
complex with coefficients in the $i$-th higher direct image
of the constant sheaf $\BQ_X$ on the \emph{Baily--Borel compactification}
of $S$. Note that the objects of \loccit \ are defined over $\BC$,
and that their behaviour under the Hecke algebra is not \emph{a priori} clear.
For $i=r \ge 1$, the image of $\Pi_S^r$ should be expected to contain
an object isomorphic to the base change of the interior motive
$\Gr_0 \Mgm(X)^e$ to $\BC$.}
Third (Corollary~\ref{4F}), the 
interior motive occurs canonically as a direct factor 
of the (Chow) motive of any smooth compactification of $X$.
Let us mention that in the case of ``non-parallel type'',
i.e., the integers $r_i$ used to construct the idempotent $e$
are not all equal to each other,
Theorem~\ref{4Main2} states that the $e$-part of the boundary motive
vanishes. In particular, the interior motive then coincides with
the $e$-part of the motive of the (open) Kuga--Sato variety.
This can be seen as a motivic explanation of \cite[Rem.~I.4.8]{F}:
``If $f$ is a modular form, but not a cusp form, then $r_1 = \ldots = r_g$.''
We then discuss the special cases $g=1$ and $g=2$. For $g=2$,
the control of the weights from Corollary~\ref{4E} turns out to be
sufficiently precise to allow for a
strenghtening of the main result of \cite{K}
(Corollaries~\ref{4I} and \ref{4K}): the \emph{special elements}
in the motivic cohomology of $X$ constructed in \loccit \ come indeed
from motivic cohomology of a (in fact, of any) smooth compactification
of $X$.  \\
 
The final Section~\ref{5} is devoted to the verification 
of the criterion from Corollary~\ref{3C} for
Hilbert--Blumenthal varieties.
First, we need to show that in this case, the boundary motive
is indeed Artin--Tate (Theorem~\ref{5A}). This is done using a 
smooth \emph{toroidal compactification}. We use \emph{co-localization}
for the boundary motive \cite{W1}, in order to reduce 
to showing the statement for the contribution of any of the strata.
The latter was identified in the general context of \emph{mixed Shimura
varieties} \cite{W2}. For Kuga--Sato families over 
Hilbert--Blumenthal varieties, \loccit \ shows 
in particular that these contributions
are indeed all Artin--Tate. We are thus reduced to 
the identification of boundary cohomology.  
The resulting formula is most certainly known to the experts
(see e.g.\ \cite{BL, Ha, Blt}). In the presence of a ``non-parallel type'', 
it actually implies the vanishing of boundary cohomology. 
In the remaining case, we employ 
the main result from \cite{BW}, to identify the weights, 
thereby completing the verification
of the criterion from Corollary~\ref{3C}.  \\

We should warn the reader that 
our constructions work \emph{a priori} with $\BQ$-coefficients.
This seems to be necessary for at least the following reasons. First,
the results on Artin--Tate motives from Section~\ref{3} 
are not known to hold before passage to $\BQ$-coefficients;
actually, it is not even clear how to formulate the integral version
of these results. 
Second, the construction of the idempotent $e$ relies on the
motivic decomposition of Abelian schemes of relative dimension greater than one
\cite{DM}. But this decomposition
necessitates the inversion of at least one prime, and 
is only known to be canonical after $\otimes \, \BQ \,$.
Finally,
our computations of the boundary cohomology of Hilbert--Blumenthal 
varieties (see Section~\ref{5}) are 
valid only after tensoring with $\BQ \,$. In fact, unless one deals
with modular curves,
very little seems to be known about the primes dividing the torsion 
of boundary cohomology with integer coefficents
(see \cite[Sect.~3.4]{Gh}).  \\

Part of this work was done while I was enjoying a 
\emph{modulation de service pour les porteurs de projets de recherche},
granted by the \emph{Universit{\'e} Paris~13}.
I wish to thank D.~Blasius, 
M.~Dimitrov, V.~Maillot, 
R.~Pink and J.~Tilouine for 
useful discussions and comments. \\

{\bf Notation and conventions}: Throughout the article, 
$k$ denotes a fixed number field. 
We denote by $Sch/k$ the category of separated schemes of finite 
type over $k$, and by $Sm/k \subset Sch/k$
the full sub-category of objects which are smooth over $k$. 
As far as motives are concerned,
the notation of this paper is that of \cite{W1,W2,W3}, which in turn follows
that of \cite{V}. We refer to \cite[Sect.~1]{W1} for a concise
review of this notation, and of the de\-fi\-nition of the triangulated 
categories $\DeffgM$ and $\DgM$ of (effective) geometrical
motives over $k$. Let $F$ be a commutative $\BQ$-algebra.
The notation $\DeffQgM$ and $\DQgM$ stands 
for the $F$-linear analogues of these triangulated categories
defined in \cite[Sect.~16.2.4
and Sect.~17.1.3]{A}. 
Similarly, let us denote by $\CHeffM$ and $\CHM$ the categories 
opposite to the categories of (effective) Chow motives, and by
$\CHeffQM$ and $\CHQM$ the pseudo-Abelian
completion of the category $\CHeffM \otimes_\BZ F$ and 
$\CHM \otimes_\BZ F$, respectively. 
Using \cite[Cor.~4.2.6]{V}, we canonically identify 
$\CHeffQM$ and $\CHQM$ with
a full additive sub-category of $\DeffQgM$ and $\DQgM \, $, respectively. \\

%%% Local Variables:
%%% mode: latex
%%% TeX-master: "head"
%%% End:

\bigskip
%\include{Sec3}

%%%%%%%%%%%%%%%%%%%%%%%%%%%%%%%%%%%%%%%%%%%%%%%%%%%%%%%%%%%%%%%%%%%%%%%
%
%  Section 3
%
%%%%%%%%%%%%%%%%%%%%%%%%%%%%%%%%%%%%%%%%%%%%%%%%%%%%%%%%%%%%%%%%%%%%%%%

\section{A criterion on the existence of the interior motive}
\label{3}

%%%%%%%%%%%%%%%%%%%%%%%%%%%%%%%%

%%%%%%%%%%%%%%%%%%%%%%%%%%%%%%%%

Fix $X \in Sm/k$, and consider the exact triangle
\[
(\ast) \quad\quad
\dMgm(X) \longto \Mgm(X) \longto \Mcgm(X) \longto \dMgm(X)[1]
\]
in $\DeffgM$. Recall from \cite[Def.~4.1~(a)]{W3} that
$c(X,X)$ contains a canonical sub-algebra $c_{1,2}(X,X)$
(of ``bi-finite correspondences'')
acting on $(\ast)$.
Denote by $\bar{c}_{1,2}(X,X)$ 
the quotient of $c_{1,2}(X,X)$ by the kernel of this action. 
Fix a finite direct product $F$ of fields of characteristic zero,
and an idempotent $e$
in $\bar{c}_{1,2}(X,X) \otimes_\BZ F$. 
Denote by $\Mgm(X)^e$, $\Mcgm(X)^e$ and $\dMgm(X)^e$ the images of $e$
on $\Mgm(X)$, $\Mcgm(X)$ and $\dMgm(X)$, respectively, considered as
objects of the category $\DeffgM_F$. 
Recall the following assumption.

\begin{Ass}[{\cite[Asp.~4.2]{W3}}] \label{3A}
The object $\dMgm(X)^e$ 
is without weights $-1$ and $0$.
\end{Ass}

Here, absence of certain weights is  
in the sense of \cite[Def.~1.10]{W3}.
In order to apply the results from \cite[Sect.~4]{W3},
allowing in particular to construct the \emph{interior motive},
one needs to verify Assumption~\ref{3A}. For this purpose,
let us consider realizations
(\cite[Sect.~2.3 and Corrigendum]{H}; see \cite[Sect.~1.5]{DG}
for a simplification of this approach). We  
shall concentrate on two
realizations (the statement from Theorem~\ref{3Main}
below then formally generalizes to
any of the other realizations ``with weights'' considered in \cite{H}):
\begin{enumerate}
\item[(i)] the Hodge theoretic realization
\[
R_\sigma : \DQgM \longto D
\]
associated to a fixed embedding $\sigma$ of the number field $k$ into
the field $\BC$ of complex numbers. Here, $D$ is the bounded derived category 
of mixed graded-polarizable $\BQ$-Hodge structures 
\cite[Def.~3.9, Lemma~3.11]{Be}, tensored with $F$,
\item[(ii)] the $\ell$-adic realization
\[
R_{\ell} : \DQgM \longto D
\]
for a prime $\ell$. Here, $D$ is the bounded ``derived category'' of
constructible $\BQ_\ell$-sheaves on $\Spec (k)$
\cite[Sect.~6]{E}, tensored with $F$.
\end{enumerate}

Choose and fix one of these two, denote it by $R$, and
recall that it is a contravariant tensor functor mapping
the pure Tate motive $\BZ(m)$ to the pure Hodge structure $\BQ(-m)$
(when $R=R_\sigma$) and to the pure $\BQ_\ell$-sheaf $\BQ_\ell(-m)$
(when $R=R_\ell$), respectively \cite[Thm.~2.3.3]{H}.
The category $D$ 
is equipped with a $t$-structure;
write $H^n$ for the cohomology functors.

\begin{Thm} \label{3Main}
Let $\alpha \le \beta$ be two integers, and $R$ one of the two realizations
considered above (Hodge theoretic or $\ell$-adic). 
Assume that $\dMgm(X)^e$ is a successive
extension of objects $M$ of $\DeffgM_F$,
each satisfying one of the following properties.
\begin{enumerate}
\item[(i)] $M$ is without weights $\alpha,\ldots,\beta$.
\item[(ii)] $M$ lies in the triangulated sub-category $\DQATM$
of $\DQgM$ of \emph{Artin--Tate motives over $k$} \cite[Def.~1.3]{W6}, and 
the cohomology object $H^n \! R(M)$ of its image $R(M)$ under $R$
is without weights $n - \beta,\ldots,n - \alpha$, 
for all $n \in \BZ$. 
\end{enumerate} 
Then $\dMgm(X)^e$ is without weights $\alpha,\ldots,\beta$.
\end{Thm}

\begin{Proof}
Apply \cite[Prop.~2.11, Thm.~3.11~(d)]{W6}.
\end{Proof}

\begin{Cor} \label{3B}
If the hypotheses of Theorem~\ref{3Main} are met with $\alpha \le -1$
and $\beta \ge 0$, then Assumption~\ref{3A} holds. 
\end{Cor}

As far as the remaining part of this article is concerned,
we shall be dealing with a situation in which the whole of
$\dMgm(X)^e$ satisfies property~(ii) from Theorem~\ref{3Main}.
It will be worthwile to spell out that property.

\begin{Cor} \label{3C}
The conclusion of Theorem~\ref{3Main} 
holds in particular if
$\dMgm(X)^e$ lies in $\DQATM$,
and if the $e$-part of the \emph{boundary cohomology} of $X$
\[
\bigl( \partial H^n (X (\BC),\BQ) \otimes_\BQ F \bigr)^e
\]
(in the Hodge theoretic setting) resp. \
\[
\bigl( \partial H^n (X_{\bar{k}},\BQ_\ell) \otimes_\BQ F \bigr)^e
\]
(in the $\ell$-adic setting) is without weights $n - \beta,\ldots,n - \alpha$, 
for all $n \in \BZ$. 
If this latter condition is fulfilled with $\alpha \le -1$
and $\beta \ge 0$, then Assumption~\ref{3A} holds.
\end{Cor}

Recall that boundary cohomology of $X$ is defined \emph{via} a compactification
$j: X \into \bX$: writing $i: \partial \! \bX \into \bX$
for the complementary immersion,
one defines $\partial H^n (\bullet)$ as cohomology of $\partial \! \bX$
with coefficients in $i^* Rj_* (\bullet)$. Thanks to proper base change,
this definition is independent of the choice of $j$,
as is the long exact cohomology sequence 
\[
\quad\quad\quad\quad\quad\quad\quad\quad\quad \ldots \longto
H^n (X (\BC),\BQ) \otimes_\BQ F \longto 
\partial H^n (X (\BC),\BQ) \otimes_\BQ F \longto 
\]
\[
H^{n+1}_c (X (\BC),\BQ) \otimes_\BQ F \longto
H^{n+1} (X (\BC),\BQ) \otimes_\BQ F \longto 
\ldots \quad\quad\quad\quad\quad\quad\quad\quad
\]
(in the Hodge theoretic setting) resp. \
\[
\quad\quad\quad\quad\quad\quad\quad\quad \ldots \longto
H^n (X_{\bar{k}},\BQ_\ell) \otimes_\BQ F \longto
\partial H^n (X_{\bar{k}},\BQ_\ell) \otimes_\BQ F \longto
\]
\[
H^{n+1}_c (X_{\bar{k}},\BQ_\ell) \otimes_\BQ F \longto
H^{n+1} (X_{\bar{k}},\BQ_\ell) \otimes_\BQ F \longto
\ldots \quad\quad\quad\quad\quad\quad\quad\quad
\]
(in the $\ell$-adic setting). \\

Note also that the algebra $\bar{c}_{1,2}(X,X)$ acts contravariantly
on the boun\-dary cohomology $\partial H^n (X (\BC),\BQ)$ 
resp.\ $\partial H^n (X_{\bar{k}},\BQ_\ell)$. \\

Corollary~\ref{3C} results from Theorem~\ref{3Main}
and the following.

\begin{Prop} \label{3D}
Fix $X \in Sm / k$, a $\BQ$-algebra $F$, and $e$ as before. Then
$H^n \! R \bigl( \dMgm(X)^e \bigr)$ is isomorphic to 
$\bigl( \partial H^n (X (\BC),\BQ) \otimes_\BQ F \bigr)^e$
(in the Hodge theoretic setting) resp. \
$\bigl( \partial H^n (X_{\bar{k}},\BQ_\ell) \otimes_\BQ F \bigr)^e$
(in the $\ell$-adic setting), for all $n$.
\end{Prop}

\begin{Proof}
It suffices to show that the image under $R$ of the
canonical morphism 
\[
\iota: \Mgm(X) \longto \Mcgm(X)
\]
can be $c_{1,2}(X,X)$-equivariantly identified with the canonical morphism
\[
R \Gamma_c(X) \longto R \Gamma(X)
\]
in the target $D$ of $R$ of classes of complexes $R \Gamma_c(X)$ and 
$R \Gamma(X)$
computing cohomology with resp.\ without support. Indeed, the 
exact triangle $(\ast)$ will then show that the $e$-part 
$\dMgm(X)^e$ of the boundary motive is mapped to a cone of
$R \Gamma_c(X)^e \to R \Gamma(X)^e$. 
It will therefore be isomorphic to the class
of a complex computing the $e$-part of boundary cohomology. 

We may assume that $X$ is of pure dimension $d$.
First note that for any fixed smooth compactification $j : X \into \bX$,
the morphism $\iota$ is the composition of the canonical
morphism
\[
\Mgm (j) : \Mgm (X) \longto \Mgm \bigl( \bX \bigr) \; ,
\]
of the inverse of the duality isomorphism, 
\[
\Mgm \bigl( \bX \bigr) \isoto \Mgm \bigl( \bX \bigr)^*(d)[2d]
\]
\cite[Thm.~4.3.7~3]{V}, of the dual of $\Mgm (j)$,
\[
\Mgm (j)^*(d)[2d] : \Mgm \bigl( \bX \bigr)^*(d)[2d]   
\longto \Mgm (X)^*(d)[2d] \; ,
\]
and of the duality isomorphism
\[
\Mgm (X)^*(d)[2d] \isoto \Mcgm (X)
\]
\cite[Thm.~4.3.7~3]{V}.
Now recall that $R$ is compatible with the tensor
structures \cite[Cor.~2.3.5, Cor.~2.3.4]{H}, and sends the Tate motive
$\BZ(1)$ to $\BQ(-1)$ resp.\ $\BQ_\ell(-1)$. 
It follows that $R$ is compatible with duality.
Furthermore, $R$ sends $\Mgm(f)$ to $f^* : R \Gamma (Z) \to R \Gamma (Y)$,
for any morphism $f: Y \to Z$ of smooth $k$-schemes \cite[pp.~6--7]{DG}.
Thus, $R$ sends $\iota$ to the composition of the duality
isomorphism
\[
R \Gamma_c (X) \isoto R \Gamma(X)^*(-d)[-2d] \; ,
\]
the dual of $j^*$, the inverse of the duality isomorphism,
\[
R \Gamma \bigl( \bX \bigr)^*(-d)[-2d] \isoto R \Gamma \bigl( \bX \bigr) \; ,
\]
and $j^*$. But this composition equals the canonical morphism
\[
R \Gamma_c(X) \longto R \Gamma(X) \; .
\]
It remains to show that the above identification is compatible with the
action of $c_{1,2}(X,X)$. Let $\FZ$ be a cycle on $X \times_k X$
belonging to $c_{1,2}(X,X)$, and denote by
${}^t \FZ$ its transpose. Both $\FZ$ and ${}^t \FZ$ are finite over
both components of $X \times X$. This is true in particular for
the first component. Therefore, both induce endomorphisms
of $\Mgm(X)$. Similarly, they induce endomorphisms of $\Mcgm(X)$.
Now the definition of the duality isomorphism
\[
\Mgm (X)^*(d)[2d] \isoto \Mcgm (X)
\]
\cite[proof of Thm.~4.3.7~3]{V} implies that under this isomorphism,
the endomorphism $\FZ$ of $\Mcgm (X)$ corresponds to
the endomorphism ${}^t \FZ^*(d)[2d]$ of $\Mgm (X)^*(d)[2d]$.
We thus identify the commutative diagram 
\[
\vcenter{\xymatrix@R-10pt{
        \Mgm (X) \ar[d]_{\FZ} \ar[r]^-{\iota} &
        \Mcgm (X) \ar[d]^{\FZ} \\
        \Mgm (X) \ar[r]^-{\iota} &
        \Mcgm (X)
\\}}
\]
with
\[
\vcenter{\xymatrix@R-10pt{
        \Mgm (X) \ar[d]_{\FZ} \ar[r]^-{\iota} &
        \Mgm (X)^*(d)[2d] \ar[d]^{{}^t \FZ^*(d)[2d]} \\
        \Mgm (X) \ar[r]^-{\iota} &
        \Mgm (X)^*(d)[2d]
\\}}
\]
By \cite[pp.~6--7]{DG}, $R$ sends $\FY$ to 
$\FY^* : R \Gamma (W) \to R \Gamma (V)$,
for any finite correspondence $\FY$ on the product $V \times_k W$
of two smooth $k$-schemes.   
It follows 
that $R$ sends the latter commutative diagram to the commutative diagram
\[
\vcenter{\xymatrix@R-10pt{
        R \Gamma (X) &
        R \Gamma (X)^*(-d)[-2d] \ar[l]  \\
        R \Gamma (X) \ar[u]^{\FZ^*} &
        R \Gamma (X)^*(-d)[-2d] \ar[l] \ar[u]_{({}^t \FZ^*)^*(-d)[-2d]}
\\}}
\] 
Now the endomorphism $({}^t \FZ^*)^*(-d)[-2d]$ of
$R \Gamma (X)^*(-d)[-2d]$ corresponds to the endomorphism $\FZ^*$ of
$R \Gamma_c (X)$. But this means precisely that our identification
of the image under $R$ of $\Mgm (X) \to \Mcgm (X)$ with 
the canonical morphism $R \Gamma_c (X) \to R \Gamma (X)$ is 
$c_{1,2}(X,X)$-equivariant.
\end{Proof}

\begin{Rem}
In the Hodge theoretic setting, the isomorphism of mixed 
Hodge structures
\[
H^n \! R \bigl( \dMgm(X) \bigr) \cong
\bigl( \partial H^n (X (\BC),\BQ) \otimes_\BQ F \bigr)
\]
from Proposition~\ref{3D}
implies an isomorphism of $F$-modules.
For quasi-projective $X$, 
this latter statement should be
compared to \cite[Lemme~3.12]{Ay}.
\end{Rem}

\forget{
\begin{Rem} \label{3E}
(a)~Consider the situation from Section~\ref{2b}, that is, fix  
a base scheme $S \in Sm/k$. 
When $X$ lies not only in $Sm/k$, but in
$Prop Sm /S$, then Theorem~\ref{2bB} provides another source
of endomorphisms of the exact triangle $(*)$, namely 
the relative Chow group $\ch^g(X \times_S X)$
(assuming $X$ to be of pure relative dimension $g$ over $S$).
In this setting, the algebra $\bar{c}_{1,2}(X,X)$ could be replaced
by $\ch^g(X \times_S X)$, and the results from this section, 
as well as those from \cite[Sect.~4]{W3} formally carry over. \\[0.1cm]
(b)~In the sequel, we shall keep the appraoch \emph{via}
$\bar{c}_{1,2}(X,X)$, because the endomorphisms ``of Hecke type''
to be considered 
will not in general come from $\ch^g(X \times_S X)$
(Example~\ref{2bL}~(b)), but still lie in the centralizer 
(in $\bar{c}_{1,2}(X,X)$) of
the idempotents we shall be interested in
(see Corollary~\ref{4Ea}). We shall
consider the subgroup
\[
\ch^g(X \times_S X)_{1,2} \subset \ch^g(X \times_S X)
\]
defined as the image of the canonical map
\[
c_S(X,X) \cap c_{1,2}(X,X) \longto \ch^g (X \times_S X) \; .
\]
\end{Rem}
}

%%% Local Variables:
%%% mode: latex
%%% TeX-master: "head"
%%% End:

\bigskip
%\include{Sec4}

%%%%%%%%%%%%%%%%%%%%%%%%%%%%%%%%%%%%%%%%%%%%%%%%%%%%%%%%%%%%%%%%%%%%%%%
%
%  Section 4
%
%%%%%%%%%%%%%%%%%%%%%%%%%%%%%%%%%%%%%%%%%%%%%%%%%%%%%%%%%%%%%%%%%%%%%%%

\section{Statement of the main results}
\label{4}

%%%%%%%%%%%%%%%%%%%%%%%%%%%%%%%%

%%%%%%%%%%%%%%%%%%%%%%%%%%%%%%%%

In order to state our main results (Theorems~\ref{4Main1}, \ref{4Main2}),
let us introduce the geometrical situation we are going to consider
from now on. The base $k$ is the field $\BQ$ of rational numbers,
$X$ is the $r$-th power of the universal Abelian scheme over a
Hilbert--Blumenthal variety of dimension $g$, and $e$ is associated  
to modular forms of weight $(r_1+2,\ldots,r_g+2)$,
for $r_1 + \ldots + r_g = r$ (see below for the precise definition). 
Theorems~\ref{4Main1} and \ref{4Main2} imply in particular
that in this context, Assumption~\ref{3A}
is satisfied as soon as $r \ge 1$: indeed,  
$\dMgm(X)$ lies in $D \! M \! AT(\BQ)_\BQ$,
and the $e$-part of the $n$-th boundary cohomology group of $X$
is without weights $n - (r-1),\ldots,n + r$, 
for all $n \in \BZ$. We then list the main consequences of 
this result (Corollaries~\ref{4E}, \ref{4Ea}, \ref{4F},
\ref{4I}, \ref{4K}), applying the theory developed in \cite[Sect.~4]{W3}.
The proofs of Theorems~\ref{4Main1} and \ref{4Main2}
will be given in Section~\ref{5}. \\
 
Fix a totally real number field $L$
of degree $g$. 
Let $I_L$ denote
the set of real embeddings of $L$.
Denote by $\ReL$
the Weil restriction from schemes over $L$ to schemes over $\BQ \,$.
The functor $\ReL$ is right adjoint to the base change 
$Z \mapsto Z_L := Z \times_\BQ L$.
Hence we have in particular a functorial adjunction morphism
$Z \to \ReL Z_L$
for any scheme $Z$ over $\BQ \, $. 
For any scheme $Y$ over $L$, 
and any subfield $F$ of $\BC$ containing the images
$\sigma(L)$ for all $\sigma \in I_L$, 
there is a canonical isomorphism
\[
(\ReL Y) \times_\BQ F \isoto \prod_{\sigma \in I_L}     
Y \times_{L,\sigma} F
\]
induced by the isomorphism
\[
L \otimes_\BQ R \isoto \prod_{\sigma \in I_L} R \; , \; 
l \otimes r \longmapsto (\sigma(l) \cdot r)_\sigma
\]
for any $F$-algebra $R$.
The composition of the adjunction (base changed by $F$)
and this isomorphism is simply the diagonal
\[
X \times_\BQ F \longto \prod_{\sigma \in I_L} X \times_\BQ F \; .
\]
By functoriality,
the determinant induces a morphism of group schemes over $\BQ$
\[
\det: \ReL \GL_{2,L} \longto \ReL \GKm \; .
\]

\begin{Def}[{\cite[Sect.~1.27]{R}}] \label{4A}
The group scheme $G$ over $\BQ$
is defined as the fibre product
\[
G := \Gm \times_{\ReL \GKm , \det} \ReL \GL_{2,L} \; .
\]
\end{Def}

In particular, we have
\[
G(\BQ) = \{ M \in \GL_2 (L) \, , \,  \det(M) \in \BQ^* \subset L^* \}  \; .
\]
Under the above isomorphism
\[
(\ReL \GL_{2,L}) \times_\BQ \BR \isoto \prod_{\sigma \in I_L}     
\GL_{2,\BR} \; ,
\]
we can identify
\[
G(\BR) \cong  
\biggl\{ (M_\sigma)_{\sigma \in I_L} \in 
\prod_{\sigma \in I_L} \GL_2 (\BR) \, , \,
\det(M_\sigma) =  \det(M_\eta) \, \forall \, 
\sigma, \eta \in I_L \biggr\}     \; .
\]
In particular, we see that $G(\BR)$ has two connected components, 
according to the sign of the determinant.
Under these identifications, the inclusion of $G(\BQ)$ into
$G(\BR)$ maps $M \in \GL_2 (L)$ to the $g$-tuple 
\[
(\sigma(M))_{\sigma \in I_L} \in 
\prod_{\sigma \in I_L} \GL_2 (\BR) \; .
\]

\begin{Def} \label{4B}
(a)~The analytic space $\CH$ is defined as
\[
\CH := \biggl\{ (\tau_\sigma)_{\sigma \in I_L} \in 
\prod_{\sigma \in I_L} (\BC - \BR) \, , \,
\sign(\imm \tau_\sigma) =  \sign(\imm \tau_\eta) \, \forall \, 
\sigma, \eta \in I_L \biggr\} \; .
\]
(b)~The action of $G(\BR)$ on $\CH$ is given by the usual componentwise action
of $\GL_2(\BR)$ on $\BC - \BR$, and the above identification 
of $G(\BR)$ with a subgroup of $\prod_{\sigma \in I_L} \GL_2 (\BR)$.
\end{Def}
   
Given that the action of $\GL_2(\BR)$ on $\BC - \BR$ is transitive
and trivial on $\Gm(\BR) \subset \GL_2(\BR)$, it is easy to see that
$G(\BR)$ acts transitively on $\CH$. Observe that this action is by
analytical automorphisms. In fact, $(G,\CH)$ are \emph{pure Shimura data}
\cite[Def.~2.1]{P}. Their \emph{reflex field} \cite[Sect.~11.1]{P}
equals $\BQ \, $. 
The center $Z(G)$ of $G$ equals
\[
\Gm \times_{\ReL \GKm , \, x \mapsto x^2} \ReL \GKm \; ,
\]
hence its neutral connected component is isogeneous to $\Gm$.
In particular, the Shimura data $(G,\CH)$ satisfy condition 
$(+)$ from \cite[Sect.~5]{W2}. \\

Let us now fix additional data:
(A)~an open compact subgroup $K$ of $G(\BA_f)$
which is neat \cite[Sect.~0.6]{P}, 
(B)~a subfield $F$ of $\BC$
containing the images of all embeddings $\sigma \in I_L$,
(C)~an integer $r \ge 0$, together with a partition
\[
\ur: \quad\quad r = \sum_{\sigma \in I_L} r_\sigma 
\]
with $g$ integers $r_\sigma \ge 0$. Equivalently,
we may see $\ur = \sum_\sigma r_\sigma \cdot \sigma$ as an element
of the free Abelian group $\BZ[I_L]$ on $I_L$. \\

These data (A)--(C) are used as follows (cmp.\ \cite[Sect.~2.2, 2.3]{K}
for the case $g = 2$).
The \emph{Shimura variety}
$S := S^K (G,\CH)$ 
is an object of $Sm / \BQ$. This is the \emph{Hilbert--Blumenthal
variety} of level $K$ associated to $L$. 
It is of dimension $g$, and
admits an interpretation as modular space of Abelian varieties
of dimension $g$ with additional structures, among which a
real multiplication by a sub-algebra of $L$
which is of rank $g$ over $\BZ$, hence of finite index in 
the ring of integers $\CO_L$ (and which depends on $K$).
In particular, there is a universal family $A$ 
of Abelian varieties over $S$. Thus,
the absolute dimension of $A$ over $\BQ$ is $2g$, it is an object of 
$Sm / \BQ \; $,
and thanks to the modular interpretation there is a canonical ring monomorphism
\[
L \longinto \End_S(A)\otimes_\BZ \BQ
\]
from $L$ to the endomorphisms of $A$ over $S$,
tensored with $\BQ \, $. Denote by
$\CHSM$ the category of \emph{Chow motives over $S$} \cite[Sect.~1.3., 1.6]{DM}. 
The decomposition of the Chow motive of $A$ over $S$
\[
h(A/S) = \bigoplus_i h_i(A/S)
\]
\cite[Thm.~3.1]{DM} being functorial, there is a map
\[
\End_S(A)\otimes_\BZ \BQ \longto 
\End_{\CHSM}\bigl( h_i(A/S) \bigr) \otimes_\BZ \BQ
\]
for all $0 \le i \le 2g$.
For $i = 1$, this map is an isomorphism 
of $\BQ$-vector spaces \cite[Prop.~2.2.1]{K}.
Hence, we get a ring monomorphism
\[
L \longinto \End_{\CHSM}\bigl( h_1(A/S) \bigr) \otimes_\BZ \BQ \; .
\]
\forget{
\begin{Lem} 
This monomorphism factors through 
\[
c_{1,2}(A,A) \otimes_\BZ \BQ \subset c(A,A) \otimes_\BZ \BQ \; .
\]
\end{Lem}

\begin{Proof}
Consider the image of $l \in L^*$ in $c(A,A) \otimes_\BZ \BQ \; $.
Multiplication by a suitable non-zero integer gives an endomorphism
of $A$.
Its graph is a cycle in $A \times_\BQ A$ which maps isomorphically
to the first component of $A \times_\BQ A$. 
Over the second component, it is necessarily
finite since $l$ is invertible. 
\end{Proof}

Thus, our monomorphism
factors to give 
\[
L \longinto c_{1,2}(A,A) \otimes_\BZ \BQ \subset 
c(A,A) \otimes_\BZ \BQ \; .
\]
} 
Its tensor product with $F$ gives
\[
L \otimes_\BQ F \longinto \End_{\CHSM}\bigl( h_1(A/S) \bigr) \otimes_\BZ F \; .
\] 
The field $F$ containing the images of all $\sigma \in I_L$, we get  
canonically
\[
L \otimes_\BQ F \isoto \prod_{\sigma \in I_L} F \; , \; 
l \otimes f \longmapsto (\sigma(l) \cdot f)_\sigma \; .
\]
In particular, there are canonical idempotents $e_\sigma$ 
in $L \otimes_\BQ F$, indexed by $I_L$: by definition, 
$e_\sigma$ is the projection to the copy of $F$ corresponding to $\sigma$. 
Let us use the same symbol $e_\sigma$ for its image in 
\[
\End_{\CHSM}\bigl( h_1(A/S) \bigr) \otimes_\BZ F \subset
\ch^g (A \times_S A) \otimes_\BZ F \; .
\] 
From our construction, the relation
\[
(\Delta) = \sum_{\sigma \in I_L} e_\sigma \in \ch^g (A \times_S A) \otimes_\BZ F
\]
is obvious. It induces a decomposition
\[
h_1(A/S) = \bigoplus_{\sigma \in I_L} h_1(A/S)^{e_\sigma} 
\]
in $\CHSQM$,
where $h_1(A/S)^{e_\sigma}$ denotes the image of 
the projector $e_\sigma$ on $h_1(A/S)$. \\

\forget{
More generally, consider the $r$-fold fibre product
$A^r := A \times_S \ldots \times_S A$ of $A$ over $S$.
\emph{Via} the partition $r = \sum_{\sigma \in I_L} r_\sigma$,
it will be identified with
\[
\rrod_{\sigma \in I_L} A^{r_\sigma} \; .
\]
The product $L^r = L \times \ldots \times L 
= \prod_{\sigma \in I_L} L^{r_\sigma}$ 
maps to $c_{1,2}(A^r,A^r) \otimes_\BZ \BQ \, $.
In 
\[
L^r \otimes_\BQ F \isoto 
\prod_{\tau \in I_L} \prod_{\sigma \in I_L} F^{r_\sigma} \; ,
\]
consider the idempotent $e_{\ur}$ corresponding to the projection
onto
\[
\prod_{\tau = \sigma} F^{r_\sigma} \; .
\]
The map $\prod_{\sigma} L^{r_\sigma} \otimes_\BQ F
\to \prod_{\sigma} F^{r_\sigma}$ induced by $e_{\ur}$ is thus given by 
\[
(l \otimes f)_\sigma \longmapsto (\sigma(l) \cdot f)_\sigma \; . 
\]
The image of $e_{\ur}$ in $c_{1,2}(A^r,A^r) \otimes_\BZ F$
will be denoted by the same symbol $e_{\ur}$,
as will its class in $\bar{c}_{1,2}(A^r,A^r) \otimes_\BZ F$. It 
commutes with all elements coming from (the commutative ring)
$L^r \otimes_\BQ F$. In particular, it commutes with
the graph $\Gamma_{[n]_{A}}$ of $[n]_{A}$, for
all integers $n$ and each of the $r$ components $A$ of $A^r$. \\

Now recall that for $0 \le i \le 2g$, we constructed
$\pi_{A,i,n} \in c_S(A,A) \otimes_\BZ \BQ$
for integers $n \ne -1, 0, 1$ (Proposition~\ref{2cD}),
which allow to decompose the
exact triangle
\[
(\ast)_A \quad\quad
\dMgm(A) \longto \Mgm(A) \longto \Mcgm(A) \longto \dMgm(A)[1] 
\]
according to the action of $\Gamma_{[n]_A}$ (Corollary~\ref{2cE}).
}
Let us now use the partition 
$\ur = \sum_\sigma r_\sigma \cdot \sigma \in \BZ[I_L]$.

\begin{Def} \label{4C}
Define ${}^{\ur} \CV \in \CHSQM$ as 
\[
{}^{\ur} \CV := 
\bigotimes_{\sigma \in I_L} \Sym^{r_\sigma} h_1(A/S)^{e_\sigma} \; .
\]
\end{Def}

The tensor product is in $\CHSQM$, and the symmetric powers are formed
with the usual convention concerning the (twist of) the natural
action of the symmetric group on a power of $A$ over $S$ (see e.g.\
\cite[p.~72]{K}). Thus, ${}^{\ur} \CV$ is a direct factor
of $h(A^r/S)$, where $A^r$ denotes the 
$r$-fold fibre product of $A$ over $S$. That is, it is associated
to an idempotent
\[
e_{\ur} \in \ch^{rg} (A^r \times_S A^r) \otimes_\BZ F \; .
\]
Let $c_S(A^r,A^r)$ denote the subgroup of $c(A^r,A^r)$ of correspondences
whose support is contained in $A^r \times_S A^r \subset A^r \times_k A^r$. 
Define 
\[
\ch^{rg} (A^r \times_S A^r)_{1,2} \subset \ch^{rg} (A^r \times_S A^r)
\]
as the image of 
\[
c_S(A^r,A^r) \cap c_{1,2}(A^r,A^r) \longto 
\ch^{rg} (A^r \times_S A^r) \; .
\]

\begin{Lem} \label{4D}
The idempotent $e_{\ur}$ lies in
\[
\ch^{rg} (A^r \times_S A^r)_{1,2} \otimes_\BZ F \subset
\ch^{rg} (A^r \times_S A^r) \otimes_\BZ F \; .
\]
\end{Lem} 

\begin{Proof}
By \cite[Prop.~3.4]{W10}, 
\[
\pi_{A,1,n} := \prod_{j \ne 1} \frac{\Gamma_{[n]_A} - n^j}{n - n^j}
\]
is a pre-image of $p_{A,1}$ in $c_S(A,A) \otimes_\BZ \BQ \,$,
for any $n \ne -1,0,1$. It visibly lies in the intersection
$(c_S(A^r,A^r) \cap c_{1,2}(A^r,A^r)) \otimes_\BZ \BQ \,$.

Similarly, $e_\sigma$ is seen to be the image
of the composition of $\pi_{A,1,n}$ and
\[
\prod_{\tau \ne \sigma} 
\frac{\Gamma_{\alpha(l)} - \tau(l)}{\sigma(l) - \tau(l)} \; ,
\]
for any $l$ generating $L$ over $\BQ$, and such that $\alpha(l)$
is a genuine endomorphism of $A$.
The graph $\Gamma_{\alpha(l)}$
is a cycle in $A \times_S A$ which maps isomorphically
to the first component of $A \times_S A$.  
Over the second component, it is necessarily
finite: indeed, the element $l$ is invertible in $L$,
hence $\alpha(l)$ in invertible in $\End_S(A)\otimes_\BZ \BQ \, $. 
Altogether, this proves that
$e_\sigma \in \ch^g (A \times_S A) \otimes_\BZ F$ comes from 
\[
\bigl( c_S(A,A) \cap c_{1,2}(A,A) \bigr) \otimes_\BZ F \; .
\]
The same is then true for the external product of the $e_\sigma$
corresponding to the direct factor  
\[
\bigotimes_{\sigma \in I_L} 
\bigl( h_1(A/S)^{e_\sigma} \bigr)^{\otimes {r_\sigma}}
\]
of $h(A^r/S)$.

In order to get a pre-image of the idempotent $e_{\ur}$,
it suffices to take a suitable average over the action of 
a suitable finite group (a product of symmetric groups).
\end{Proof}

According to \cite[Cor.~2.12]{W10}, the 
idempotent $e_{\ur}$ thus maps to an idempotent in 
$\bar{c}_{1,2}(A^r,A^r) \otimes_\BZ F$. It will
be denoted by the same symbol $e_{\ur}$. By functoriality
\cite[Thm.~2.2~(a)]{W10}, the relative Chow motive
\[
{}^{\ur} \CV = h(A^r/S)^{e_{\ur}}
\]
gives rise to an exact triangle
\[
\dMgm(A^r)^{e_{\ur}} \longto \Mgm(A^r)^{e_{\ur}} \longto \Mcgm(A^r)^{e_{\ur}} 
\longto \dMgm(A^r)^{e_{\ur}}[1] 
\]
in $\DeffQgM$.
Here are our main results.

\begin{Thm} \label{4Main1}
The boundary motive $\dMgm(A^r)$ lies in the
triangulated sub-category $D \! M \! DT(\BQ)_\BQ$ of 
$DM_{gm}(\BQ)_\BQ$ of \emph{Dirichlet--Tate motives 
over $\BQ \,$} \cite[Def.~3.5~(b)]{W6}. Its direct factor
$\dMgm(A^r)^{e_{\ur}}$ 
is without weights 
\[
-r, -(r-1), \ldots, r-1 \; .
\]
In particular, Assumption~\ref{3A} holds 
for $\dMgm(A^r)^{e_{\ur}}$ whenever $r \ge 1$.
\end{Thm}

\begin{Thm} \label{4Main2}
Assume that there are $\tau, \sigma \in I_L$ such that $r_\tau \ne r_\sigma$
(hence $g \ge 2$ and $r \ge 1$).
Then $\dMgm(A^r)^{e_{\ur}} = 0$, and 
$\Mgm(A^r)^{e_{\ur}} \cong \Mcgm(A^r)^{e_{\ur}}$
are effective Chow motives.
\end{Thm}

Theorems~\ref{4Main1} and \ref{4Main2}
will be proved in Section~\ref{5}.  
Let us give their main corollaries,
assuming that $r \ge 1$.
First,
we fix a weight filtration 
\[ 
C_{\le -(r+1)} \longto \dMgm(A^r)^{e_{\ur}} \longto 
C_{\ge r} \longto C_{\le -(r+1)}[1]
\]
avoiding weights $-r, \ldots, r-1$ \cite[Def.~1.6]{W3}. The category $D \! M \! DT(\BQ)_\BQ$
being pseudo-Abelian \cite[Cor.~2.6]{W6}, the motives
$C_{\le -(r+1)}$ and $C_{\ge r}$ are Dirichlet--Tate motives over $\BQ$
of weights $\le -(r+1)$ and $\ge r$, respectively.
Furthermore, $C_{\le -(r+1)} = 0 = C_{\ge r}$
under the hypothesis of Theorem~\ref{4Main2}.

\begin{Cor}[{\cite[Thm.~4.3]{W3}}]  \label{4E}
Assume $r \ge 1$. \\[0.1cm]
(a)~The motive $\Mgm(A^r)^{e_{\ur}}$ is without weights $-r, \ldots, -1$,
and the motive $\Mcgm(A^r)^{e_{\ur}}$ is without weights $1, \ldots, r$. 
The Chow motives $\Gr_0 \Mgm(A^r)^{e_{\ur}}$ 
and $\Gr_0 \Mcgm(A^r)^{e_{\ur}}$ \cite[Prop.~2.2]{W3} are defined,
and they carry a natural action of 
\[
GCen_{\bar{c}_{1,2}(A^r,A^r)}(e_{\ur}) :=
\big{ \{ } z \in \bar{c}_{1,2}(A^r,A^r) \otimes_\BZ F \; , \;
ze_{\ur} = e_{\ur}ze_{\ur} \big{ \} } \; .
\]
(b)~There are canonical exact triangles 
\[
C_{\le -(r+1)} \longto \Mgm(A^r)^{e_{\ur}} \stackrel{\pi_0}{\longto}
\Gr_0 \Mgm(A^r)^{e_{\ur}} \longto C_{\le -(r+1)}[1]
\]
and
\[
C_{\ge r} \longto \Gr_0 \Mcgm(A^r)^{e_{\ur}} \stackrel{i_0}{\longto}
\Mcgm(A^r)^{e_{\ur}} \longto C_{\ge r}[1] \; ,
\]
which are stable under the natural action
of $GCen_{\bar{c}_{1,2}(A^r,A^r)}(e_{\ur})$. \\[0.1cm]
(c)~There is a canonical isomorphism 
$\Gr_0 \Mgm(A^r)^{e_{\ur}} \isoto \Gr_0 \Mcgm(A^r)^{e_{\ur}}$ in $\CHeffM_F$. 
As a morphism, it is uniquely determined by the property
of making the diagram
\[
\vcenter{\xymatrix@R-10pt{
        \Mgm(A^r)^{e_{\ur}} \ar[r]^-{u} \ar[d]_{\pi_0} &
        \Mcgm(A^r)^{e_{\ur}} \\
        \Gr_0 \Mgm(A^r)^{e_{\ur}} \ar[r] &
        \Gr_0 \Mcgm(A^r)^{e_{\ur}} \ar[u]_{i_0}
\\}}
\]
commute; in particular, it is 
$GCen_{\bar{c}_{1,2}(A^r,A^r)}(e_{\ur})$-equivariant. \\[0.1cm]
(d)~Let $N \in \CHM_F$ be a Chow motive. Then $\pi_0$ and $i_0$ induce
isomorphisms
\[
\Hom_{\CHM_F} \bigl( \Gr_0 \Mgm(A^r)^{e_{\ur}} , N \bigr) \isoto
\Hom_{\DgM_F} \bigl( \Mgm(A^r)^{e_{\ur}} , N \bigr)
\]
and
\[
\Hom_{\CHM_F} \bigl( N , \Gr_0 \Mcgm(A^r)^{e_{\ur}} \bigr) \isoto
\Hom_{\DgM_F} \bigl( N , \Mcgm(A^r)^{e_{\ur}} \bigr) \; .
\]
(e)~Let $\Mgm(A^r)^{e_{\ur}} \to N \to \Mcgm(A^r)^{e_{\ur}}$ 
be a factorization of $u$ through a Chow motive
$N \in \CHM_F$. 
Then $\Gr_0 \Mgm(A^r)^{e_{\ur}} = \Gr_0 \Mcgm(A^r)^{e_{\ur}}$ is canonically
a direct factor of $N$, with a canonical direct complement.
\end{Cor} 

Henceforth, we identify $\Gr_0 \Mgm(A^r)^{e_{\ur}}$ and 
$\Gr_0 \Mcgm(A^r)^{e_{\ur}}$
\emph{via} the canonical isomorphism of Corollary~\ref{4E}~(c).
Note that under the hypothesis of Theorem~\ref{4Main2},
we have 
\[
\Mgm(A^r)^{e_{\ur}} = \Gr_0 \Mgm(A^r)^{e_{\ur}} \quad \text{and} \quad
\Gr_0 \Mcgm(A^r)^{e_{\ur}} = \Mcgm(A^r)^{e_{\ur}} \; ,
\]
and the isomorphism of Corollary~\ref{4E}~(c)
coincides with that of Theorem~\ref{4Main2}.
The equivariance statements from
Corollary~\ref{4E}~(a)--(c) apply in particular
to cycles coming from the \emph{Hecke algebra} associated to
the Shimura variety $S$.
More precisely, we have the following statement.

\begin{Cor} \label{4Ea}
Assume $r \ge 1$. Then $\Gr_0 \Mgm(A^r)^{e_{\ur}}$
carries a natural action of the Hecke algebra $R(K,G(\BA_f))$
associated to the neat open compact subgroup $K$ of $G(\BA_f)$.
More precisely, any $x \in G(\BA_f)$ defines a cycle denoted $KxK$ in 
$c_{1,2}(A^r,A^r)$, whose class in $\bar{c}_{1,2}(A^r,A^r)$ belongs to
the centralizer 
\[
Cen_{\bar{c}_{1,2}(A^r,A^r)}(e_{\ur}) :=
\big{ \{ } z \in \bar{c}_{1,2}(A^r,A^r) \otimes_\BZ F \; , \;
ze_{\ur} = e_{\ur}z \big{ \} } 
\subset GCen_{\bar{c}_{1,2}(A^r,A^r)}(e_{\ur})
\]
of $e_{\ur}$.
\end{Cor}

\begin{Proof}
Fix $x \in G(\BA_f)$. Recall that our base scheme $S$ equals  
the Hilbert--Blumenthal variety $S^K(G,\CH)$. It is the target
of two finite \'etale morphisms $g_1,g_2: U \to S$, where
$U$ denotes the Hilbert--Blumenthal variety $S^{K \cap x^{-1}Kx}(G,\CH)$.
Using the notation of \cite[Sect.~3.4]{P}, 
the morphism $g_1$ equals $[\ \cdot 1]$, and the morphism $g_2$
equals $[\ \cdot x^{-1}]$. Note that by the very definition of the
$\BQ$-rational structure of $S$ and $U$ (e.g.\ \cite[Def.~11.5]{P}),
both $g_1$ and $g_2$ are indeed defined over $\BQ \, $.

Recall that $A$ is the universal Abelian scheme over $S$;
denote by $A_1, A_2$ its base changes to $U$ \emph{via} $g_1$ and $g_2$,
respectively. To the data $K$ and $x$, the following are
canonically associated: a third Abelian scheme $B$ over $U$
admitting real multiplication,
and isogenies $f_1: B \to A_1$ and $f_2: B \to A_2$
compatible with the real multiplications.  
By definition, the cycle $KxK$ is then equal to the direct image under
$g_1 \times_k g_2$ of the composition
\[
\Gamma_{f_2^r} \circ {}^t \Gamma_{f_1^r} 
\in c_U (A_1^r,A_2^r) \cap c_{1,2} (A_1^r,A_2^r)
\]
($=$ ``pull-back \emph{via} $f_1^r$ followed by push-out \emph{via} $f_2^r$'').

In order to show that the class of $KxK$  
in $\bar{c}_{1,2}(A^r,A^r)$
commutes with $e_{\ur} \, $, note first that
$\varphi := \Gamma_{f_2^r} \circ {}^t \Gamma_{f_1^r}$
defines a morphism of relative Chow motives over $U$,
\[
\varphi: h(A_1^r / U) = g_1^* ( h(A^r/S) ) \longto 
g_2^* ( h(A^r/S) ) = h(A_2^r / U) \; .
\] 
Then,
since both $f_1$ and $f_2$ are isogenies, this morphism
is compatible 
with the external products of the idempotents $p_{A_i,1}$
\cite[Thm.~3.1, Prop.~3.3]{DM}. Since $f_1$ and $f_2$ also respect the real
multiplication, the morphism $\varphi$ is also compatible with the cycle classes
$g_i^*(e_{\ur}) = e_{\ur,i} 
\in \ch^{rg} (A^r_i \times_U A^r_i) \otimes_\BZ F$,
i.e., we have the relation
\[
\varphi \circ g_1^*(e_{\ur}) = g_2^*(e_{\ur}) \circ \varphi
\]
of morphisms of relative Chow motives over $U$.
We are thus in the situation of \cite[Ex.~2.16~(d)]{W10}.
Now observe that using the notation from \loccit, 
the effect of the cycle $KxK$ on the exact triangle
\[
\dMgm(A^r) \longto \Mgm(A^r) 
\stackrel{u}{\longto} \Mcgm(A^r) 
\longto \dMgm(A^r)[1] 
\]
coincides with $\varphi(g_1,g_2)$ \cite[Ex.~2.16~(e)]{W10}.
Thanks to the relation 
$\varphi \circ g_1^*(e_{\ur}) = g_2^*(e_{\ur}) \circ \varphi$,
we thus get that the class of $KxK$ in $\bar{c}_{1,2}(A^r,A^r)$
belongs indeed to 
$Cen_{\bar{c}_{1,2}(A^r,A^r)}(e_{\ur})$.
\end{Proof}

\begin{Cor}[{\cite[Cor.~4.6]{W3}}]  \label{4F} 
Assume $r \ge 1$, and
let $\widetilde{A^r}$ be any smooth compactification of $A^r$. Then
$\Gr_0 \Mgm(A^r)^{e_{\ur}}$ is canoni\-cally
a direct factor of the Chow motive
$\Mgm(\widetilde{A^r})$, with a canonical direct complement.
\end{Cor}

Furthermore, \cite[Thm.~4.7, Thm.~4.8]{W3} on the Hodge theoretic
and $\ell$-adic realizations \cite[Cor.~2.3.5, Cor.~2.3.4 and Corrigendum]{H}
apply, and tell us in particular that $\Gr_0 \Mgm(A^r)^{e_{\ur}}$
is mapped to the part of interior cohomology of $A^r$
fixed by $e_{\ur}$. In particular, the $L$-function of the 
Chow motive $\Gr_0 \Mgm(A^r)^{e_{\ur}}$
is computed \emph{via} (the $e_{\ur}$-part of) interior cohomology of $A^r$.

\begin{Def}[{\cite[Def.~4.9]{W3}}] \label{4G}
Let $r \ge 1$. 
We call $\Gr_0 \Mgm(A^r)^{e_{\ur}}$ the \emph{$e_{\ur}$-part of the
interior motive of $A^r$}.
\end{Def}

\begin{Rem} \label{4H}
By \cite[Thm.~4.14]{W3}, control of the reduction of \emph{some}
compactification of $A^r$ implies control of certain properties
of the $\ell$-adic realization of $\Gr_0 \Mgm(A^r)^{e_{\ur}}$.
To the best of the author's knowledge, the sharpest result known 
about reduction of compactifications of $A^r$ is \cite[Thm.~6.4]{DimT}.
It concerns the case when $K \subset G(\BA_f)$ is of type
$\Gamma_1$, and states that there exist then smooth compactifications
of $A^r$ having good reduction at each prime number $p$
dividing neither the level $N$ of $K$ nor the absolute discriminant
$d$ of $L$. 
\cite[Thm.~4.14]{W3} then yields the following conclusions: 
(a)~for all primes $p$ not dividing $Nd$, the $p$-adic
realization of $\Gr_0 \Mgm(A^r)^{e_{\ur}}$ is crystalline,
(b)~if furthermore $p \ne \ell$, then
the $\ell$-adic realization of $\Gr_0 \Mgm(A^r)^{e_{\ur}}$
is unramified. Note that
given the identification of the $\ell$-adic realization of 
$\Gr_0 \Mgm(A^r)^{e_{\ur}}$ with intersection cohomology,
conclusions~(a) ane (b) are already contained
in \cite[Sect.~7]{DimT}. 
\forget{
Furthermore, \cite[proof of Thm.~2.3]{Dim}
shows how to transfer statements~(a) and (b) 
to Hilbert--Blumenthal varieties associated to
more general open compact neat subgroups $K$ of $G(\BA_f)$. 
}
\end{Rem}

\begin{Rem} \label{4Ha}
(a)~If all $r_\sigma$ are strictly positive (hence $r \ge g$) 
then Saper's vanishing theorem on (ordinary) cohomology 
\cite[Thm.~5]{Sa} implies that the realizations of $\Gr_0 \Mgm(A^r)^{e_{\ur}}$
are concentrated in the single cohomological degree $r+g$.  
In particular, we expect the following relation to the Chow--K\"unneth
decompositions constructed in \cite[Thm.~2.4]{GHM}. 
The base change from $\BQ$ to $\BC$ of $\Gr_0 \Mgm(A^r)^{e_{\ur}}$
should map monomorphically to the $(r+g)$-th Chow--K\"unneth component
of the motive (over $\BC$) of any toroidal com\-pactification of $A^r$. \\[0.1cm]
(b)~In general, consider the (relative) Chow--K\"unneth projectors
$\Pi_S^i$, $i = 0,\ldots,2rg$ of \cite[Thm.~I]{GHM1} modelling intersection cohomology of $S(\BC)$
with coefficients in the $i$-th higher direct image of the constant sheaf $\BQ_{A^r}$. 
For $i=r \ge 1$, we expect the image of $\Pi_S^r$ to contain a copy
of the base change to $\BC$ of the interior motive
$\Gr_0 \Mgm(X)^e$.
\end{Rem}

Let us discuss special cases. First, for $g=1$, we have $L = \BQ \,$,
$S$ is a smooth modular curve, and $A$ the universal family of
elliptic curves over $S$ (see our Introduction). As coefficient field~(B),
we may choose $F = \BQ \, $. The partition~(C) amounts to fixing an integer
$r \ge 0$, and 
${}^{\ur} \CV = \Sym^r h_1(A/S)$. In this setting,
the formal implications of Theorem~\ref{4Main1} are discussed in
\cite[Rem.~4.17]{W3}; note that our idempotent $e_{\ur}$
coincides with the idempotent denoted $e$ in \loccit. 
Indeed, the additional action of torsion entering the definition of $e$
is known (and ea\-sily shown)
to be trivial on the relative motive $h(A/S)$. 
In particular \cite[Rem.~4.17~(b)]{W3}, we get an alternative 
construction of the Grothendieck
motive $M(f)$ associated to a normalized newform $f$ of weight $r+2$
\cite{Sch} as a direct factor of the Grothendieck motive
underlying $\Gr_0 \Mgm(A^r)^{e_{\ur}}$. \\

Now let $g = 2$. Here, we have $[L : \BQ] = 2$,
$S$ is a smooth Hilbert--Blumenthal surface, and $A$ the universal family of
Abelian surfaces over $S$. As coefficient field~(B),
we may choose $F$ equal to the image of $L$ under any of its two 
real embeddings $\tau, \sigma \in I_L$. 
The partition~(C) amounts to fixing two integers
$r_\tau, r_\sigma \ge 0$, whose sum is denoted $r$. Then 
\[
{}^{\ur} \CV = 
\Sym^{r_\tau} h_1(A/S)^{e_\tau} \otimes \Sym^{r_\sigma} h_1(A/S)^{e_\sigma} \; .
\]
For any object $M$ of $DM_{gm}^{eff}(\BQ)_F $, define motivic cohomology 
\[
H_{\CM}^p \bigl( M,F(q) \bigr) :=  
\Hom_{DM_{gm}^{eff}(\BQ)_F} \bigl( M,\BZ(q)[p] \bigr) \; .
\]
When $M = \Mgm(Y)$ for a scheme $Y \in Sm / \BQ$, 
this gives motivic cohomology
$H_{\CM}^p \bigl( Y,\BZ(q) \bigr)$ of $Y$,
tensored with $F$. Now observe that the
relative Chow motive ${}^{\ur} \CV$ coincides with the object denoted
$\CV^{r_\tau,r_\sigma}_K$ in \cite[Def.~2.3.1]{K}.
From now on, assume that $r_\tau \ge r_\sigma \ge 1$
(hence $r \ge 2$). The main result of \loccit \
gives the construction of a sub-space
\[
\CK(r_\tau,r_\sigma,n) \subset H^{r+3}_{\CM} \bigl( {}^{\ur} \CV , F(n) \bigr)
= \Hom_{DM_{gm}^{eff}(\BQ)_F} \bigl( \Mgm(A^r)^{e_{\ur}} , \BZ(n)[r+3] \bigr)
\]
for all integers $n$ between $r_\tau + 2$ and $r+2 = r_\tau + r_\sigma + 2$
\cite[Thm.~5.2.4]{K}, and establishes a weak version of Beilinson's
conjecture for Asai $L$-functions \cite[Thm.~5.2.4~(b)]{K}.
As already mentioned in \cite[p.~62, Rem.~5.2.5~(a)]{K},
one of the shortcomings of this result is that $\CK(r_\tau,r_\sigma,n)$
is not shown to come from motivic cohomology of a 
smooth compactification of $A^r$. It is reasonable to expect this to be true;
one of the indications being
\cite[Thm.~5.2.4~(a)]{K} that
the Hodge theoretic realization of $\CK(r_\tau,r_\sigma,n)$ lands in
\[
H^{r+3}_{! \FH} \bigl( {}^{\ur} \CV_{/ \BR} , \BR(n) \bigr) \subset
H^{r+3}_{\FH} \bigl( {}^{\ur} \CV_{/ \BR} , \BR(n) \bigr) \; ,
\] 
which by definition \cite[Def.~(2.4.1)]{K} is the sub-space
of absolute Hodge cohomology of ${}^{\ur} \CV_{/ \BR}$ given by the image of
absolute Hodge cohomology of any smooth compactification of $A^r$. 
Our main results allow to give a significantly more precise statement.

\begin{Cor} \label{4I}
Assume that $r_\tau \ge r_\sigma + 1 (\ge 2)$
or that $r_\tau + 2 \le n \le r + 1$.
Then the map on the level of motivic cohomology induced by the
morphism $\pi_0: \Mgm(A^r)^{e_{\ur}} \to \Gr_0 \Mgm(A^r)^{e_{\ur}}$,
\[
\pi_0^*: 
\Hom_{DM_{gm}(\BQ)_F} \bigl( \Gr_0 \Mgm(A^r)^{e_{\ur}} , \BZ(n)[r+3] \bigr)
\longto H^{r+3}_{\CM} \bigl( {}^{\ur} \CV , F(n) \bigr)
\]
is an isomorphism.
\end{Cor}  

In particular, under the hypotheses of the corollary, 
$\CK(r_\tau,r_\sigma,n)$ can be considered as a sub-space of
$\Hom_{DM_{gm}(\BQ)_F} ( \Gr_0 \Mgm(A^r)^{e_{\ur}} , \BZ(n)[r+3] )$,
and hence (Corollary~\ref{4F}) of motivic cohomology of any 
smooth compactification of $A^r$.

\begin{Cor} \label{4K}
Assume that $r_\tau = r_\sigma (\ge 1)$, hence $r = 2r_\tau$.
Then the image of the map on the level of motivic cohomology induced by
$\pi_0$,
\[
\pi_0^* \bigl( \Hom_{DM_{gm}(\BQ)_F} 
\bigl( \Gr_0 \Mgm(A^r)^{e_{\ur}} , \BZ(r+2)[r+3] \bigr) \bigr)
\subset H^{r+3}_{\CM} \bigl( {}^{\ur} \CV , F(r+2) \bigr)
\]
contains the sub-space $\CK(r_\tau,r_\sigma,r+2)$.
\end{Cor}  

In particular, under the hypotheses of the corollary, 
$\CK(r_\tau,r_\sigma,r+2)$ comes from a sub-space of
$\Hom_{DM_{gm}(\BQ)_F} ( \Gr_0 \Mgm(A^r)^{e_{\ur}} , \BZ(r+2)[r+3] )$.

\begin{Rem}
This settles the problem raised in \cite[Rem.~5.2.5~(a)]{K}.
At least two other points remain, in order to get a proof
of Beilinson's full conjecture:
first \cite[Rem.~5.2.5~(c)]{K}, the elements in $\CK(r_\tau,r_\sigma,n)$
should be \emph{integral} (with respect to suitable models over $\Spec \BZ$),
second, the space of integral elements in motivic cohomology
should be equal to $\CK(r_\tau,r_\sigma,n)$.
We have nothing to say about these two points. 
\end{Rem}

\medskip

\begin{Proofof}{Corollaries~\ref{4I} and \ref{4K}}
Recall the exact triangle
\[
C_{\le -(r+1)} \longto \Mgm(A^r)^{e_{\ur}} \stackrel{\pi_0}{\longto}
\Gr_0 \Mgm(A^r)^{e_{\ur}} \longto C_{\le -(r+1)}[1]
\]
from Corollary~\ref{4E}~(b). Here,
$C_{\le -(r+1)}$ is a Dirichlet--Tate motive over $\BQ$
of weights $\le -(r+1)$.
Theorem~\ref{4Main2} tells us that
\[
C_{\le -(r+1)} = 0 
\]
when $r_\tau \ge r_\sigma + 1$. In this case, the morphism $\pi_0$
is therefore itself an isomorphism, and thus induces an isomorphism
on the level of motivic cohomology.

In general, the kernel of 
\[
\pi_0^*: 
\Hom_{DM_{gm}(\BQ)_F} \bigl( \Gr_0 \Mgm(A^r)^{e_{\ur}} , \BZ(n)[r+3] \bigr)
\longto H^{r+3}_{\CM} \bigl( {}^{\ur} \CV , F(n) \bigr)
\]
is a quotient of
\[
\Hom_{DM_{gm}(\BQ)_F} \bigl( C_{\le -(r+1)}[1] , \BZ(n)[r+3] \bigr) =
\Hom_{DM_{gm}(\BQ)_F} \bigl( C_{\le -(r+1)} , \BZ(n)[r+2] \bigr) 
\]
and the co-kernel a sub-space of
\[
\Hom_{DM_{gm}(\BQ)_F} \bigl( C_{\le -(r+1)} , \BZ(n)[r+3] \bigr) \; .
\]
Now $\BZ(n)$ is pure of weight $-2n$. 
When $n \le r+1$, then the weights of $\BZ(n)[r+2]$ 
and of $\BZ(n)[r+3]$ are at least equal to $-2(r+1) + r+2 = -r$. 
Since $C_{\le -(r+1)}$ is of weights at most $-(r+1)$, 
orthogonality \cite[Def.~1.1.1~(iii)]{Bo} of the motivic weight structure
\cite[Sect.~6.5 and 6.6]{Bo}
implies that both
$\Hom_{DM_{gm}(\BQ)_F} ( C_{\le -(r+1)} , \BZ(n)[r+2] )$ and  
$\Hom_{DM_{gm}(\BQ)_F} ( C_{\le -(r+1)} , \BZ(n)[r+3] )$
are zero. In this case, the map $\pi_0^*$ is therefore again an isomorphism. 

In the sequel, let us therefore assume that $r_\tau = r_\sigma$,
and that $n = r+2$.
As above, the co-kernel injects into
\[
\Hom_{DM_{gm}(\BQ)_F} \bigl( C_{\le -(r+1)} , \BZ(n)[r+3] \bigr) \; .
\]
Let us first show that on this latter space, the map induced by the
Hodge theoretic realization is injective. 
Observe that the 
motive $\BZ(n)[r+3]$ is pure of weight $-(r+1)$.
This is the highest weight possibly occurring in $C_{\le -(r+1)}$.
Shifting by $r+1$ therefore reduces us to show the following:
for any two Dirichlet--Tate motives $M$ and $N$, with
$M \in D \! M \! DT(\BQ)_{F, w \le 0}$ and $N \in D \! M \! DT(\BQ)_{F, w = 0}$,
the map induced on
\[
\Hom_{DM_{gm}(\BQ)_F} ( M , N ) 
\]
by the Hodge theoretic realization is injective.
Choose an exact triangle
\[
M_{\le -2} \longto M \longto M_{-1,0} \longto M_{\le -2}[1] \; ,
\]
with $M_{\le -2} \in D \! M \! DT(\BQ)_{F, w \le -2}$ and 
\[
M_{-1,0} \in D \! M \! DT(\BQ)_{F, w \ge -1} 
\cap D \! M \! DT(\BQ)_{F, w \le 0} \; .
\] 
By \cite[Cor.~2.9]{W6}, the object $M_{-1,0}$ is a direct sum
of two Dirichlet--Tate motives $M_{-1} \oplus \Gr_0 M$,
the first being pure of weight $-1$ and the second pure of weight $0$. 
Orthogonality formally implies
that the two morphisms
$M \to M_{-1,0}$ and $M_{-1,0} \to \Gr_0 M$ induce isomorphisms
\[
\Hom_{DM_{gm}(\BQ)_F} ( M_{-1,0} , N ) \isoto \Hom_{DM_{gm}(\BQ)_F} ( M , N ) 
\]
and
\[
\Hom_{DM_{gm}(\BQ)_F} ( \Gr_0 M , N ) 
\isoto \Hom_{DM_{gm}(\BQ)_F} ( M_{-1,0} , N ) \; .
\]
The Hodge theoretic realization $R$ maps our data to an exact triangle
\[
R(M_{\le -2}) \longleftarrow R(M) \longleftarrow R(M_{-1,0}) 
\longleftarrow R(M_{\le -2})[-1] 
\]
in the bounded derived category $D$ of mixed graded-polarizable Hodge structures
(recall that $R$ is contravariant),
and a direct sum decomposition
\[
R(M_{-1}) \oplus R(\Gr_0 M) = R(M_{-1,0})\; .
\]
According to \cite[Rem.~3.13~(a)]{W6},
the functor $R$ respects the weight structures. This means that
$R(M_{\le -2})$ has weights at least $2$, that $R(M_{-1})$ is pure of 
weight $1$ and that $R(\Gr_0 M)$ and $R(N)$ are pure of weight $0$.
As above, orthogonality (for the category $D$) yields
formally that 
\[
\Hom_D ( K , R(M_{-1,0}) ) \isoto \Hom_D ( K , R(M) )    
\]
and
\[
\Hom_D ( K , R(\Gr_0 M) ) \isoto \Hom_D ( K, R(M_{-1,0}) ) 
\]
for any object $K$ of the heart $D_{w = 0}$, hence in particular
for $K = R(N)$. Altogether, we are thus reduced to showing injectivity of
\[
R: \Hom_{DM_{gm}(\BQ)_F} ( M , N ) 
\longto \Hom_D ( R(N) , R(M) )
\]
under the additional assumption that $M$ belongs to the heart 
$D \! M \! DT(\BQ)_{F, w = 0}$, too. In other words, we must show
faithfulness of the restriction of $R$ to $D \! M \! DT(\BQ)_{F, w = 0}$.
But this follows easily from the explicit description of
$D \! M \! DT(\BQ)_{F, w = 0}$ given in \cite[Thm.~2.5~(c)]{W6},
and the (obvious) faithfulness of $R$ on the category $\QDM$ 
from \cite[Def.~3.5~(a)]{W6}.

In order to finish the proof, it remains to show that the image
of the space $\CK(r_\tau,r_\sigma,n)$ in 
\[
\Hom_{DM_{gm}(\BQ)_F} \bigl( C_{\le -(r+1)} , \BZ(n)[r+3] \bigr) 
\]
is mapped to zero under $R$.
Choose a smooth compactification $\widetilde{A^r}$ of $A^r$. 
By Corollary~\ref{4E}~(d), the morphism
$j: \Mgm(A^r)^{e_{\ur}} \to \Mgm(A^r) \to \Mgm(\widetilde{A^r})$
factors through $\Gr_0 \Mgm(A^r)^{e_{\ur}}$.
It follows that the exact triangle 
\[
C_{\le -(r+1)} \longto \Mgm(A^r)^{e_{\ur}} \longto
\Gr_0 \Mgm(A^r)^{e_{\ur}} \longto C_{\le -(r+1)}[1]
\]
maps to an exact triangle of the form
\[
C' \stackrel{i}{\longto} \Mgm(A^r)^{e_{\ur}} \stackrel{j}{\longto}
\Mgm(\widetilde{A^r}) \longto C'[1] \; .
\]
By \cite[Thm.~5.2.4~(a)]{K}, the sub-space 
\[
\CK(r_\tau,r_\sigma,n) \subset 
\Hom_{DM_{gm}^{eff}(\BQ)_F} \bigl( \Mgm(A^r)^{e_{\ur}} , \BZ(n)[r+3] \bigr)
\]
vanishes under the composition of
\[
i^* : 
\Hom_{DM_{gm}^{eff}(\BQ)_F} \bigl( \Mgm(A^r)^{e_{\ur}} , \BZ(n)[r+3] \bigr)
\longto \Hom_{DM_{gm}^{eff}(\BQ)_F} \bigl( C' , \BZ(n)[r+3] \bigr)
\]
and of $R$. \emph{A fortiori}, its image in
\[
\Hom_{DM_{gm}(\BQ)_F} \bigl( C_{\le -(r+1)} , \BZ(n)[r+3] \bigr) 
\]
vanishes under $R$.
\end{Proofof}

\begin{Rem} \label{4L}
The picture for $r_\tau = r_\sigma$ and $n = r+2$ remains incomplete:
it is clearly desirable to identify a \emph{canonical}
pre-image of $\CK(r_\tau,r_\sigma,r+2)$ in
\[
\Hom_{DM_{gm}(\BQ)_F} 
\bigl( \Gr_0 \Mgm(A^r)^{e_{\ur}} , \BZ(r+2)[r+3] \bigr) \; .
\]
In order to achieve this, one needs to go into the construction
of the elements from \cite{K}. The vital ingredient is Beilinson's
\emph{Eisenstein symbol} \cite{Be1}, which needs to be re-interpreted
in the context of the category $DM_{gm}(\BQ)_F$.  
We plan to treat this elsewhere.
\end{Rem}

\begin{Rem}
The case $r=0$ is not covered by our main results. It concerns
the motive and the motive with compact support of the base scheme $S$.
In this situation, the best replacement 
of the interior motive is the (hypothetical) \emph{intersection motive}
$M^{!*}(S)$ of $S$ with respect to its Baily--Borel compactification. 
Note that for $g \le 2$, its existence is established:
for modular curves ($g = 1$), it equals the motive of \emph{the} smooth
compactification $S^*$ of $S$. For $g=2$, we refer to \cite{CM,W}. 
\end{Rem}

\bigskip
%\include{Sec5}

%%%%%%%%%%%%%%%%%%%%%%%%%%%%%%%%%%%%%%%%%%%%%%%%%%%%%%%%%%%%%%%%%%%%%%%
%
%  Section 5
%
%%%%%%%%%%%%%%%%%%%%%%%%%%%%%%%%%%%%%%%%%%%%%%%%%%%%%%%%%%%%%%%%%%%%%%%

\section{Proof of the main results}
\label{5}

%%%%%%%%%%%%%%%%%%%%%%%%%%%%%%%%

%%%%%%%%%%%%%%%%%%%%%%%%%%%%%%%%

We keep the notation of the preceding section. In order to prove
Theorems~\ref{4Main1} and \ref{4Main2}, our idea is to apply the
criterion from Section~\ref{3}. Let us start by fixing some notation.

\begin{Def} 
Denote by $V$ the standard
two-dimensional representation of $\GL_{2,L}$ over $L$.
\end{Def}

Thus, $\ReL V$ is a $2g$-dimensional representation of 
$\ReL \GL_{2,L}$, and therefore, of $G$. 

\begin{Thm} \label{5A}
For any integer $r \ge 0$,
the boundary motive $\dMgm(A^r)$ lies in the
triangulated sub-category $D \! M \! DT(\BQ)_\BQ$ of 
$DM_{gm}(\BQ)_\BQ$ of Dirichlet--Tate motives 
over $\BQ \,$.
\end{Thm}

\begin{Proof}
The variety $A^r$ is a \emph{mixed Shimura variety} over 
$S = S^K (G,\CH)$. More precisely, the
representation $\ReL V$ of $G$ is easily seen to be of Hodge type
$\{ (-1,0) , (0,-1) \}$ in the sense of \cite[Sect.~2.16]{P}.
The same statement is then true for the $r$-th power $\ReL V^r$
of $\ReL V$. By \cite[Prop.~2.17]{P}, this allows for the construction
of the \emph{unipotent extension} $(P^r,\FX^r)$ of $(G,\CH)$
by $\ReL V^r$. The reader wishing an explicit description of
$(P^r,\FX^r)$ is referred to \cite[Sect.~1.1]{K}, where the case 
$g=2$ is treated. The description from \loccit \ generalizes
easily to arbitrary $g$.

The pair $(P^r,\FX^r)$ constitute \emph{mixed
Shimura data} \cite[Def.~2.1]{P}. By construction, they come endowed with
a morphism $\pi^r: (P^r,\FX^r) \to (G,\CH)$ of Shimura data,
identifying $(G,\CH)$ with the pure Shimura data
underlying $(P^r,\FX^r)$. In particular, $(P^r,\FX^r)$ also satisfy
condition $(+)$ from \cite[Sect.~5]{W2}.

Now there is an open compact neat subgroup $K^r$ of $P^r(\BA_f)$,
whose image under $\pi^r$ equals $K$, and
such that $A^r$ is identified with the \emph{mixed Shimura variety}
$S^{K^r} (P^r,\FX^r)$ \cite[Sect.~3.22, Thm.~11.18 and 11.16]{P}.
Furthermore, the morphism $\pi^r$ of Shimura data induces a morphism
$S^{K^r} (P^r,\FX^r) \to S^K (G,\CH)$, which is identified with the
structure morphism of $A^r$.

In order to obtain control on
the boundary motive of $A^r$, we fix a smooth
\emph{toroidal compactification} $\widetilde{A^r}$. It is associated to a 
\emph{$K^r$-admissible complete smooth cone 
decomposition} $\FS$, i.e., a collection of subsets of
\[
\CC (P^r,\FX^r) \times P^r (\BA_f)
\]
satisfying the axioms of \cite[Sect.~6.4]{P}. Here,
$\CC (P^r,\FX^r)$
denotes the \emph{co\-ni\-cal complex} associated to $(P^r,\FX^r)$
\cite[Sect.~4.24]{P}. 

We refer to \cite[9.27, 9.28]{P} for criteria sufficient to guarantee
the existence of the associated compactification
$\widetilde{A^r} := S^{K^r} (P^r , \FX^r , \FS)$. 
It comes equipped with a natural (finite) 
stratification into locally closed strata.
The unique open stratum is $A^r$. Any stratum 
$\widetilde{A^r_\sigma}$ different from the generic
one is associated to a \emph{rational boundary component} $(P_1,\FX_1)$
of $(P^r,\FX^r)$ \cite[Sect.~4.11]{P} which is \emph{proper},
i.e., unequal to $(P^r,\FX^r)$. 

First, co-localization for the boundary motive \cite[Cor.~3.5]{W1}
tells us that $\dMgm(A^r)$ is a successive extension of 
(shifts of) objects of the form
\[
\Mgm (\widetilde{A^r_\sigma}, i_\sigma^! \, j_! \, \BZ) \; .
\]
Here, $j$ denotes the open immersion of $A^r$ into $\widetilde{A^r}$,
$i_\sigma$ runs through the immersions of 
the strata $\widetilde{A^r_\sigma}$ different from $A^r$ into $\widetilde{A^r}$,
and $\Mgm (\widetilde{A^r_\sigma}, i_\sigma^! \, j_! \, \BZ)$
is the motive of $\widetilde{A^r_\sigma}$ 
with coefficients in $i_\sigma^! \, j_! \, \BZ$
defined in \cite[Def.~3.1]{W1}. 

Next, by \cite[Thm.~6.1]{W2}, there is an isomorphism
\[
\Mgm (\widetilde{A^r_\sigma}, i_\sigma^! \, j_! \, \BZ) \isoto
\Hom (\BZ (\sigma),\Mgm (S^{K_1}(P_1,\FX_1))) [\dim \sigma]  \; .
\]
Recall \cite[p.~971]{W2} that the \emph{group of orientations}
$\BZ (\sigma)$ is (non-canonically) isomorphic to $\BZ$,
hence 
\[
\Hom (\BZ (\sigma),\Mgm (S^{K_1}(P_1,\FX_1))) \cong 
\Mgm (S^{K_1}(P_1,\FX_1))  \; .
\]
$S^{K_1}(P_1,\FX_1)$ is a Shimura variety associated to the data
$(P_1,\FX_1)$ and an open compact neat subgroup $K_1$ of $P_1(\BA_f)$. 
In order to show our claim, we are thus reduced to showing that
$\Mgm (S^{K_1})$ is an object of $D \! M \! DT(\BQ)_\BQ$, 
for any Shimura variety
$S^{K_1} = S^{K_1}(P_1,\FX_1)$ associated to a proper rational
boundary component $(P_1,\FX_1)$ of $(P^r,\FX^r)$, and any
open compact neat subgroup $K_1$ of $P_1(\BA_f)$.

Given that $P^r$ is a unipotent extension of $G$,
the pure Shimura data underlying $(P_1,\FX_1)$ coincides
with the pure Shimura data underlying some proper rational
boundary component $(G_1,\CH_1)$ of $(G,\CH)$. By definition
\cite[Sect.~4.11]{P}, the group $G_1$ is associated
to an \emph{admissible $\BQ$-parabolic subgroup} $Q$ of $G$
\cite[Def.~4.5]{P}.
It is not difficult to see that the inverse image under the immersion
of $G$ into $\ReL \GL_{2,L}$ induces a bijection on the 
sets of $\BQ$-parabolic subgroups. Under this bijection, the group
$Q$ corresponds necessarily 
to a group of the form $\ReL B$, for some Borel subgroup $B$
of $\GL_{2,L}$. Equivalently, $Q$ is the stabilizer in $G$ of a
subspace of $\ReL V$ of the form $\ReL V'$, for a one-dimensional
$L$-subspace $V'$ of $V$. A computation ana\-logous to the one
from \cite[Ex.~4.25]{P} for $g=1$ shows that the pure Shimura data 
underlying $(G_1,\CH_1)$ equal the data $(\BG_{m,\BQ},\CH_0)$
from \cite[Ex.~2.8]{P}. Altogether, we see that $(P_1,\FX_1)$
is a unipotent extension of $(\BG_{m,\BQ},\CH_0)$.

We are ready to conclude.
As follows directly
from the definition of the canonical model
(cmp.~\cite[Sect.~11.3, 11.4]{P}),
the pure Shimura variety $S^{\pi^r(K_1)}(\BG_{m,\BQ},\CH_0)$ 
underlying $S^{K_1}$ equals the spectrum of a cyclotomic field $C$
over $\BQ \, $. By \cite[Prop.~11.14]{P}, the variety $S^{K_1}$
is isomorphic to a power of the multiplicative group
over $C$. In particular, its motive lies in $D \! M \! DT(\BQ)_\BQ \, $. 
\end{Proof}

\begin{Rem} \label{5Aa}
As the proof shows, Theorem~\ref{5A} admits a version ``before tensoring
with $\BQ \,$''. That is, the boundary motive $\dMgm(A^r)$ lies in the
triangulated sub-category $D \! M \! DT(\BQ)$
of $DM_{gm}(\BQ)$ generated by Tate twists and motives $\Mgm(\Spec k)$,
for Abelian finite field extensions $k$ of $\BQ \,$.
\end{Rem}

In order to apply the results from Section~\ref{3}, we need to analyze
the Hodge structure on the $e_{\ur}$-part of the boundary cohomology of $A^r$,
\[
\bigl( \partial H^n \bigl( A^r (\BC),\BQ \bigr) 
                   \otimes_\BQ F \bigr)^{e_{\ur}} \; ,
\]
for all integers $n$. Recall that there is a canonical isomorphism
\[
(\ReL V) \otimes_\BQ F \isoto \bigoplus_{\sigma \in I_L} V_\sigma \; ,
\]
where we set $V_\sigma := V \otimes_{L,\sigma} F$.
In fact, this is an isomorphism of representations of $G$ over $F$.

\begin{Def} 
Denote by $V_{\ur}$ the representation
\[
V_{\ur} := 
\bigotimes_{\sigma \in I_L} \Sym^{r_\sigma} V_\sigma^{\vee} \; .
\]
of $G$ over $F$.
\end{Def}

The tensor product is over $F$, and $V_\sigma^{\vee}$ is the 
contragredient representation of $V_\sigma$.
Recall (e.g.\ \cite[Thm.~2.2]{LNM}) the definition of
the \emph{canonical construction} functor $\mu$ from the category
of finite-dimensional algebraic representations of $G$ to the category
of admissible graded-polarizable variations of Hodge structure on $S$.

\begin{Prop} \label{5B}
There is a canonical isomorphism of Hodge structures
\[
\bigl( \partial H^n \bigl( A^r (\BC),\BQ \bigr) \otimes_\BQ F \bigr)^{e_{\ur}} 
\isoto
\partial H^{n-r} \bigl( S (\BC), \mu(V_{\ur}) \bigr) 
\]
for all integers $n$.
\end{Prop}

\begin{Proof}
The central observation is that the functor $\mu$ maps $\ReL V^\vee$
to the first higher direct image of $\BQ(0)$ under the structure morphism
of the Abelian variety $A$ \cite[remark following Lemma~2.5]{LNM}.
The rest of the argument is purely formal; cmp.\ \cite[proof of Prop.~2.3.3]{K}
for the case $g=2$.
\end{Proof}

Thus, we need to control $\partial H^{n-r} ( S (\BC), \mu(V_{\ur}) )$. 
As the reader may expect, we use the \emph{Baily--Borel compactification} 
$S^*$ of $S$. The complement of $S$ consists
of finitely many cusps;
the boundary cohomology of $S (\BC)$ therefore coincides with the direct sum
over the cusps of the \emph{degeneration} of the coefficients to the 
boundary of $S^*$. 

\medskip

\begin{Proofof}{Theorem~\ref{4Main2}}
By assumption, there are $\tau, \sigma \in I_L$ 
such that $r_\tau \ne r_\sigma$. This implies
\cite[Lemme~2.2.8]{BL} that the boundary cohomology of $\mu(V_{\ur})$
vanishes. Proposition~\ref{3D} tells us that the 
Hodge theoretic realization of
$\dMgm(A^r)^{e_{\ur}}$ is zero.
But $\dMgm(A^r)^{e_{\ur}}$ is Artin--Tate (Theorem~\ref{5A}
and \cite[Cor~2.6]{W6}). Given that the realization is conservative
on Artin--Tate motives \cite[Cor.~3.10~(a)]{W6},
this means that $\dMgm(A^r)^{e_{\ur}}$ is itself zero.
\end{Proofof}

\medskip

\begin{Proofof}{Theorem~\ref{4Main1}}
The first claim is Theorem~\ref{5A}.
Given Theorem~\ref{4Main2}, we may assume that all $r_\sigma$ are equal,
\[
r_\sigma = s \ge 0 \; \forall \, \sigma \in I_L \; ,
\] 
say. Thus,
\[
V_{\ur} = 
\bigotimes_{\sigma \in I_L} \Sym^s V_\sigma^{\vee} 
\]
and $r = g \cdot s$. Note that this representation 
descends to $\BQ \,$, and that it occurs as a direct factor of the
representation $\Sym^r \ReL V^{\vee}$. Fix a cusp $y$ of $S^* (\BC)$, 
and denote by
$j$ the open immersion of $S$ into $S^*$. 
We need to compute the weights occurring in 
\[
R^{n-r} j_* \bigl(  \mu(V_{\ur}) \bigr)_y \; .
\]
First \cite[bottom of p.~386]{BL}, the cup product
\[
R^0 j_* \bigl(  \mu(V_{\ur}) \bigr)_y \otimes_\BQ R^{n-r} j_* ( \BQ(0) )_y 
\longto R^{n-r} j_* \bigl(  \mu(V_{\ur}) \bigr)_y
\]
is an isomorphism in degrees $0 \le n-r \le g-1$. 
Next \cite[Thm.~1.3.4, Cor.~1.3.7]{BL}, in the same range of indices,
the map induced by the cup product
\[
\Lambda^{n-r} R^1 j_* ( \BQ(0) )_y \longto R^{n-r} j_* ( \BQ(0) )_y
\] 
is an isomorphism. We shall show:
\begin{enumerate}
\item[(1)] $R^0 j_* \bigl(  \mu(V_{\ur}) \bigr)_y \cong \BQ(0)$
as Hodge structures,
\item[(2)] $R^1 j_* ( \BQ(0) )_y \cong \BQ(0)^{g-1}$
as Hodge structures, when $g \ge 2$.
\end{enumerate}

Admitting these claims for the moment, we see from (1) and (2) that
the Hodge structure
$R^{n-r} j_* ( \mu(V_{\ur}) )_y$ is pure of weight $0$ for 
$0 \le n-r \le g-1$,
i.e., $r \le n \le r+g-1$. In particular, 
$R^{n-r} j_* ( \mu(V_{\ur}) )_y$ is without weights
$n-(r-1),\ldots,n+r$ whenever $n \le r+g-1$. 
To deal with the complementary range of indices $n \ge r+g$,
recall that the Hodge structures 
\[
R^m j_* \bigl(  \mu(V_{\ur}) \bigr)_y \quad \text{and} \quad
R^{2g-1-m} j_* \bigl(  \mu(V_{\ur}) \bigr)_y(r+g)
\]
are dual to each other, for all integers $m$. Indeed, the representation
$V_{\ur}$ is pure of weight $r$ in the sense of \cite[Sect.~1.11]{P},
therefore \cite[Prop.~1.12]{P} it can be $G$-equivariantly identified
with its own contragredient, twisted by $\BQ(-r)$, i.e.,
\[
\mu(V_{\ur}) \cong \mu(V_{\ur})^\vee(-r)
\]
as variations of Hodge structures. Given the definition
of the Verdier dual $\BD$ in the category of algebraic Hodge modules
\cite[Prop.~2.6]{Sai}, we have
\[
\mu(V_{\ur}) \cong \BD_S (\mu(V_{\ur})) (-(r+g))[-2g] \; .
\]
Denoting by $i$ the immersion of $y$ into $S^*$, we have
\[
i^* \circ Rj_* \circ \BD_S = \BD_y \circ i^! \circ j_! \; ;
\]
furthermore, $i^! \circ j_! = i^* \circ Rj_*[1]$ \cite[formulae~(4.3.5) and (4.4.1)]{Sai}.
This implies that
\[
i^* Rj_* \mu(V_{\ur}) \cong \BD_y i^* Rj_* \mu(V_{\ur}) (-(r+g))[-2g+1] \; .
\]
Therefore, we see that $R^{n-r} j_* ( \mu(V_{\ur}) )_y$ 
is pure of weight $2(r+g)$ whenever 
$g \le n-r \le 2g-1$,
i.e., $r+g \le n \le r+2g-1$. In particular, 
$R^{n-r} j_* ( \mu(V_{\ur}) )_y$ is without weights
$n-(r-1),\ldots,n+r$ in any case. Given Theorem~\ref{5A},
we thus have verified the hypotheses of Corollary~\ref{3C}. 

It remains to show claims (1) and (2). We shall
use the main result from 
\cite{BW} on degeneration in the Baily--Borel compactification  
of variations in the image of $\mu$. The cusp $y$ belongs to 
one of the strata associated to a rational boundary component
$(P_1,\FX_1)$ of $(G,\CH)$, where $P_1$ is contained
as a normal subgroup in one of the admissible $\BQ$-parabolic
subgroups $Q$ of $G$ \cite[Sect.~4.11]{P}. The latter being the stabilizer
in $G$ of a subspace of $\ReL V$ of the form $\ReL V'$, for a one-dimensional
$L$-subspace $V'$ of $V$, we see that the situation is conjugate under
an element of $G(\BQ)$ to the one associated to the standard Borel subgroup.
Since claims~(1)
and (2) are invariant under isomorphisms, we may therefore assume that
we work in this setting. 

It is identical to the one considered in \cite[proof of Prop.~3.2]{Blt}
(with the same notation). Applying \cite[Thm.~2.9]{BW}, we see that
\[
R^0 j_* \bigl(  \mu(V_{\ur}) \bigr)_y \cong 
H^0 \bigl( \bar{H}_C , H^0(W_1, \Res_Q^G V_{\ur}) \bigr) \; ,
\]
while 
\[
R^1 j_* ( \BQ(0) )_y \cong
H^0 \bigl( \bar{H}_C , H^1(W_1, \BQ(0)) \bigr) \oplus
H^1 \bigl( \bar{H}_C , H^0(W_1, \BQ(0)) \bigr) \; .
\]
Here, $W_1$ denotes the unipotent radical of $Q$, and $\bar{H}_C$
is free Abelian of rank $g-1$ (cmp.\ \cite[Sect.~3.2]{Blt}).  
Note that there is a shift by the codimension of $y$ in $S^*$
(which equals $g$)
in the formula of \cite[Thm.~2.9]{BW}, due to the normalization
of the inclusion of the category of variations into the derived
category of algebraic Hodge modules used in \loccit .
In order to evalute the first of the above expressions, 
one proceeds dually to
\cite[proof of Prop.~3.2]{Blt}, to show:
\begin{enumerate}
\item[(3)] $H^0(W_1, \Res_Q^G V_{\ur})$ is one-dimensional,
\item[(4)] the actions of $\bar{H}_C$ and of $P_1/W_1$ on 
$H^0(W_1, \Res_Q^G V_{\ur})$ are both trivial.
\end{enumerate}
Given that the action of $P_1/W_1$ determines the Hodge structure,
this shows (1).
As for (2), note first that 
\[
H^0 \bigl( \bar{H}_C , H^1(W_1, \BQ(0)) \bigr) = 0
\]
when $g \ge 2$ \cite[bottom of p.~386]{BL}. Hence
\[
R^1 j_* ( \BQ(0) )_y \cong
H^1 \bigl( \bar{H}_C , H^0(W_1, \BQ(0)) \bigr) 
\]
in this case.
Given that $H^0(W_1, \BQ(0)) = \BQ(0)$, and that  
the action of $\bar{H}_C$ on $H^0(W_1, \BQ(0))$
is trivial (use claim~(4) for $s=0$), we have indeed
\[
R^1 j_* ( \BQ(0) )_y \cong
\Hom \bigl( \bar{H}_C , \BQ \bigr)(0) \; .
\]
\end{Proofof}

\bigskip

%%%%%%%%%%%%%%%%%%%%%%%%%%%%%%%%%%%%%%%%%%%%%%%%%%%%%%%%%%%%%%%%%%%%%%%
%
%  Bibliography
%
%%%%%%%%%%%%%%%%%%%%%%%%%%%%%%%%%%%%%%%%%%%%%%%%%%%%%%%%%%%%%%%%%%%%%%%

\end{document}